\documentclass[11pt,reqno]{amsart}
\usepackage{amssymb,verbatim,enumerate,ifthen}
\usepackage[mathscr]{eucal}
\usepackage[utf8]{inputenc}
\usepackage[T1]{fontenc}
\def\N{\mathbb{N}}
\def\R{\mathbb{R}}
\def\Q{\mathbb{Q}}
\def\Z{\mathbb{Z}}

\def\RX{\overline{\R}}

\def\M{\mathscr{M}}
\def\AM{\mathbb{A}}
\def\MM{\mathbb{M}}
\def\D{\mathscr{D}}

\def\AC{\mathscr{C}}

\def\lip{\mathop{\mbox{\rm Lip}}\nolimits}
\def\id{\mathop{\mbox{\rm id}}\nolimits}

\def\uplus{\,\dot{+}\,}
\def\lplus{\,\text{\d{\ensuremath+}}\,}

\def\inw{\mathop{\mbox{\rm inw}}\nolimits}

\newtheorem{theorem}{Theorem}[section]
\newtheorem*{theorem*}{Theorem}
\long\def\Thm#1#2{\ifthenelse{\equal{#1}{*}}{\begin{theorem*}#2\end{theorem*}}
             {\begin{theorem}\label{T#1}#2\end{theorem}}}
\newtheorem{Atheorem}{Theorem}

\def\thm#1{Theorem~\ref{T#1}}

\newtheorem{proposition}[theorem]{Proposition}
\newtheorem*{proposition*}{Proposition}
\long\def\Prp#1#2{\ifthenelse{\equal{#1}{*}}{\begin{proposition*}#2\end{proposition*}}
             {\begin{proposition}\label{P#1}#2\end{proposition}}}
\def\prp#1{Proposition~\ref{P#1}}

\newtheorem{corollary}[theorem]{Corollary}
\newtheorem*{corollary*}{Corollary}
\long\def\Cor#1#2{\ifthenelse{\equal{#1}{*}}{\begin{corollary*}#2\end{corollary*}}
             {\begin{corollary}\label{C#1}#2\end{corollary}}}
\def\cor#1{Corollary~\ref{C#1}}

\newtheorem{lemma}[theorem]{Lemma}
\newtheorem*{lemma*}{Lemma}
\long\def\Lem#1#2{\ifthenelse{\equal{#1}{*}}{\begin{lemma*}#2\end{lemma*}}
             {\begin{lemma}\label{L#1}#2\end{lemma}}}
\def\lem#1{Lemma~\ref{L#1}}

\theoremstyle{definition}
\newtheorem{definition}[theorem]{Definition}
\newtheorem*{definition*}{Definition}
\long\def\Defn#1#2{\ifthenelse{\equal{#1}{*}}{\begin{definition*}\rm #2\end{definition*}}
             {\begin{definition}\label{D#1}\rm #2\end{definition}}}

\newtheorem{remark}[theorem]{Remark}
\newtheorem*{remark*}{Remark}
\long\def\Rem#1#2{\ifthenelse{\equal{#1}{*}}{\begin{remark*}\rm #2\end{remark*}}
             {\begin{remark}\label{R#1}\rm #2\end{remark}}}

\newtheorem{example}[theorem]{Example}
\newtheorem*{example*}{Example}
\long\def\Exa#1#2{\ifthenelse{\equal{#1}{*}}{\begin{example*}\rm #2\end{example*}}
             {\begin{example}\label{Ex#1}\rm #2\end{example}}}

\def\eq#1{{\rm(\ref{E#1})}}
\def\Eq#1#2{\ifthenelse{\equal{#1}{*}}
  {\begin{equation*}\begin{aligned}[]#2\end{aligned}\end{equation*}}
  {\begin{equation}\begin{aligned}[]\label{E#1}#2\end{aligned}\end{equation}}}

\begin{document}
\begin{flushright}
\end{flushright}
\vspace{5mm}

\date{\today}

\title{Implications between generalized convexity properties of real functions}

\author[T. Kiss]{Tibor Kiss}
\author[Zs. P\'ales]{Zsolt P\'ales}
\address{Institute of Mathematics, University of Debrecen,
H-4010 Debrecen, Pf.\ 12, Hungary}
\email{\{kiss.tibor,pales\}@science.unideb.hu}

\subjclass[2000]{Primary 39B52, Secondary 46C99}
\keywords{convexity with respect to a mean; descendant of means; fixed point theorems; divided differences}

\thanks{This research has been supported by the Hungarian
Scientific Research Fund (OTKA) Grant K111651}

\begin{abstract}
Motivated by the well-known implications among $t$-convexity properties of real functions, analogous relations 
among the upper and lower $M$-convexity properties of real functions are established. More precisely, 
having an $n$-tuple $(M_1,\dots,M_n)$ of continuous two-variable means, the notion of the descendant of these 
means (which is also an $n$-tuple $(N_1,\dots,N_n)$ of two-variable means) is introduced. In particular, 
when all the means $M_i$ are weighted arithmetic, then the components of their descendants are also weighted 
arithmetic means. More general statements are obtained in terms of the generalized quasi-arithmetic or 
Matkowski means. The main results then state that if a function $f$ is $M_i$-convex for all 
$i\in\{1,\dots,n\}$, then it is also $N_i$-convex for all $i\in\{1,\dots,n\}$. Several consequences are 
discussed.
\end{abstract}

\maketitle

\section{Introduction}

In the theory of convex functions the notion of $t$-convexity plays an important role. For $t\in\,]0,1[\,$ a real 
function $f:I\to\R$ (where $I$ is a nonempty real interval) is termed \textit{$t$-convex} (cf.\ Kuhn 
\cite{Kuh84}), Nikodem--P\'ales \cite{NikPal03}) if, for all $x,y\in I$, the inequality
\Eq{*}{
  f(tx+(1-t)y)\leq tf(x)+(1-t)f(y)
}
hold. The $\frac12$-convex functions are usually called \textit{Jensen convex}. If a function is 
$t$-convex for all $t\in\,]0,1[\,$ then it is simply called \textit{convex}. Among the many implications 
related to $t$-convexity properties we mention the following ones:
\begin{enumerate}\itemsep=2mm
 \item If $f$ is Jensen convex then it is $\Q$-convex, i.e., $t$-convex for all $t\in[0,1]\cap\Q$ 
(Kuczma \cite{Kuc85});
 \item If $f$ is $t$-convex for some $t\in\,]0,1[\,$, then it is Jensen convex (Dar\'oczy--P\'ales 
\cite{DarPal87});
 \item If $f$ is $t$-convex for some $t\in\,]0,1[\,$, then, by a result of Kuhn \cite{Kuh84}, there exists a 
subfield $K$ of $\R$ such that
 \Eq{*}{
  \{s\in\,]0,1[\,\mid \mbox{$f$ is $s$-convex}\}=\,]0,1[\,\cap\,K.
 }
\end{enumerate}
For more general results related to higher-order convexity notions refer to the paper by Gil\'anyi and P\'ales 
\cite{GilPal08}. 

We recall now the notion of \textit{second-order divided difference} defined for $f:I\to\R$ and pairwise distinct 
elements $x,y,z$ of $I$ by
\Eq{*}{
 [x,y,z;f]
 :=\frac{f(x)}{(y-x)(z-x)}+\frac{f(y)}{(x-y)(z-y)}+\frac{f(z)}{(x-z)(y-z)}.
}
In terms of this concept, the $t$-convex functions have the following easy-to-see characterization: \textit{A 
function $f:I\to\R$ is $t$-convex if and only if, for all $x,y\in I$ with $x\neq y$, we have 
$[x,tx+(1-t)y,y;f]\geq0$.}

In the paper Nikodem--P\'ales \cite{NikPal03}, $t$-convex functions were also characterized by the 
nonnegativity of a certain second-order derivative which is analogous to the standard characterization of 
twice differentiable convex functions. In this paper, an inequality related to the second-order divided 
differences was also established which turned out to be a key tool for the proofs of the main results therein.

\Prp{chi}{\textup{(Chain Inequality)}
Let $I\subseteq\R$ be an interval and $f:I\to\R$. Then, for all $n\in\N$, $x_{0}<x_{1}<\dots<x_{n+1}$ 
in $I$, and for all $i\in\{1,\dots, n\}$, the following inequalities hold:
\Eq{*}{
\min\limits_{1\,\leq\,j\,\leq\,n}[x_{j-1},x_{j},x_{j+1};f]
\leq [x_{0},x_{i},x_{n+1};f]
\leq\max\limits_{1\,\leq\,j\,\leq\,n}[x_{j-1},x_{j},x_{j+1};f].
}}

To demonstrate the use of this inequality, we show that $t$-convexity implies Jensen convexity for every real 
function $f:I\to\R$. Assume that $f:I\to\R$ is $t$-convex for some $t\in\,]0,\frac12[\,$ and let $x,y\in I$
with $x<y$ be arbitrary points. Set 
\Eq{*}{
  x_0:=x,\qquad x_1:=tx+(1-t)\frac{x+y}{2},\qquad 
  x_2:=\frac{x+y}{2},\qquad x_3:=t\frac{x+y}{2}+(1-t)y,\qquad
  x_4:=y.
}
Then
\Eq{*}{
  x_1=tx_0+(1-t)x_2,\qquad x_2=tx_3+(1-t)x_1,\qquad x_3=tx_2+(1-t)x_4,
}
whence, by the $t$-convexity of $f$, we have
\Eq{*}{
  [x_0,x_1,x_2;f]\geq0,\qquad [x_1,x_2,x_3;f]\geq0,\qquad [x_2,x_3,x_4;f]\geq0.
}
In view of the Chain Inequality, this implies that $[x_0,x_2,x_4;f]\geq0$ also holds, which is equivalent 
to the Jensen convexity of $f$.

The Jensen convexity property of a function is equivalent to the restricted condition
\Eq{*}{
  f\big(\tfrac12x+\tfrac12y\big)\leq\tfrac12f(x)+\tfrac12f(y) \qquad (x,y\in I,\, x<y). 
}
On the other hand, for $t\in\,]0,1[\,\setminus\{\frac12\}$, the $t$-convexity property is equivalent to the 
condition
\Eq{*}{
  \left\{\begin{array}{c}
  f(tx+(1-t)y)\leq tf(x)+(1-t)f(y) \\[2mm] f((1-t)x+ty)\leq (1-t)f(x)+tf(y)
  \end{array}\right. \qquad (x,y\in I,\, x<y),
}
that is, the $t$-convexity property can be expressed in terms of two inequalities over the triangle 
$\{(x,y)\in I^2\mid x<y\}$. It turns out that these two inequalities are not consequences of each other for 
every $t\,]0,1[\,\setminus\{\frac12\}$. In 2014, for every transcendental number $t$, Lewicki and Olbry\'s 
\cite{LewOlb14} constructed a function $f:I\to\R$ such that 
\Eq{*}{
  \left\{\begin{array}{c}
  f(tx+(1-t)y) < tf(x)+(1-t)f(y) \\[2mm] f((1-t)x+ty) > (1-t)f(x)+tf(y)
  \end{array}\right. \qquad (x,y\in I,\, x<y).
}
It is, however, unknown if these two inequalities are equivalent to each other for rational, or more 
generally, for algebraic $t$. (Moreover, the particular case $t=\frac13$ also has not been answered yet.) 

In this paper, for a given two-variable mean $M:\{(x,y)\in I^2\mid x\leq y\}\to\R$, we consider the class 
of functions $f:I\to\R$ satisfying the inequality
\Eq{*}{
  f(M(x,y))\leq\frac{y-M(x,y)}{y-x}f(x)+\frac{M(x,y)-x}{y-x}f(y) 
}
for all $x,y\in I$ with $x<y$. Such functions will be called $M$-convex. In this terminology, the 
$t$-convexity of a function $f:I\to\R$ is equivalent to its convexity with respect to the means $\AM_t$ and 
$\AM_{1-t}$, where, for
$s\in[0,1]$, the \textit{weighted arithmetic mean} $\AM_s:\{(x,y)\in \R^2\mid x\leq y\}\to\R$ defined by  
\Eq{*}{
\AM_s(x,y)=sx+(1-s)y.
}
Observe that, for $x<y$ and $0<s<t<1$, we have
\Eq{*}{
  x=\min(x,y)=\AM_1(x,y)<\AM_t(x,y)<\AM_s(x,y)<\AM_0(x,y)=\max(x,y)=y. 
}
Motivated by the above-described implications among $t$-convexity properties, we are going to establish 
analogous relations among the $\AM_t$-convexity properties of real functions. More generally, we will 
introduce and investigate the notions of upper and lower $M$-convexity for extended real valued functions. The 
main results of the paper then establish several implications between these convexity properties. 

More precisely, having an $n$-tuple $(M_1,\dots,M_n)$ of continuous means, we introduce the notion of 
the descendant of these means which is also an $n$-tuple $(N_1,\dots,N_n)$ of means. In several cases, we 
explicitly construct the descendant of a given $n$ tuple of means. In particular, 
when all the means $M_i$ are weighted arithmetic then the components of their descendants are also weighted 
arithmetic means. More general statements are also obtained in terms of the generalized quasi-arithmetic or 
Matkowski means. In our main results we then prove that if a function $f$ is $M_i$-convex for all 
$i\in\{1,\dots,n\}$, then it is also $N_i$-convex for all $i\in\{1,\dots,n\}$.

\section{Notations and terminology}

If $n,m\in\Z$ then set $\{k\in\Z\mid n\leq k\text{ and }k\leq m\}$ will be denoted by $\{n,\dots,m\}$. According 
to this convention $\{n,\dots,m\}=\emptyset$ if $m<n$ and $\{n,\dots,m\}$ is the singleton $\{n\}$ if $n=m$.

Given a subset $S\subseteq\R$ and $n\in\N$, we denote the set of increasingly and strictly increasingly 
ordered $n$-tuples of $S$ by $S_\leq^n$ and $S_<^n$, i.e.,
\Eq{*}{
  S_\leq^n:=\{(t_1,\dots,t_n)\in S^n\mid t_1\leq\cdots\leq t_n\}
  \qquad\mbox{and}\qquad
  S_<^n:=\{(t_1,\dots,t_n)\in S^n\mid t_1<\cdots<t_n\},
}
respectively.

A function $M:S_\leq^2\to\R$ is called a \emph{two-variable mean on $S$} and a \emph{two-variable strict mean 
on $S$} if 
\Eq{*}{
  x \leq M(x,y) \leq y \qquad \big((x,y)\in S_\leq^2\big) 
  \qquad\mbox{and}\qquad  
  x < M(x,y) < y \qquad \big((x,y)\in S_<^2\big),
}
respectively. We note that, two-variable means are usually defined on the Cartesian product $S^2$, however, 
in our approach the values of means on $S_>^2:=S^2\setminus S_\leq^2$ are irrelevant. Obviously, if 
$T\subseteq S$, then the restriction $M|_{T_\leq^2}$ is also a mean on $T$.

In the subsequent sections $I$ always denotes a nonempty interval of $\R$.

The most important class of two-variable means that appears in the consequences of our results is the class 
of generalized quasi-arithmetic means introduced by J.\ Matkowski \cite{Mat10b} in 2010: We say that a 
function $M:I_\leq^2\to\R$ is a \textit{generalized quasi-arithmetic mean} or a \textit{Matkowski mean} if 
there exist continuous, strictly increasing functions $f,g:I\to\R$ such that
\Eq{*}{
M(x,y)=\MM_{f,\,g}(x,y):=(f+g)^{-1}(f(x)+g(y)) \qquad \big((x,y)\in I_\leq^2\big).
}
Under the conditions of this definition, it is obvious that $\MM_{f,\,g}$ is a continuous strict 
mean on $I$ which is also strictly increasing in each of its variables. 

If $s\in\,]0,1[$ and $f:I\to\R$ is a continuous strictly increasing function then the Matkowski mean 
$\MM_{s f,\,(1-s)f}$ is called a \textit{weighted quasi-arithmetic mean}. We can see that
\Eq{*}{
\MM_{sf,\,(1-s)f}(x,y)=f^{-1}(sf(x)+(1-s)f(y))\qquad \big((x,y)\in I_\leq^2\big).
}
For $s=1/2$ this function is termed a \textit{(symmetric) quasi-arithmetic mean}.
Finally, observe that the mean $\MM_{s\,\mathrm{id},\,(1-s)\,\mathrm{id}}$ equals the weighted 
arithmetic mean $\AM_s$. 

\section{Auxiliary results}

\Thm{pd}{
For $n\in\N$ and for the vectors $u=(u_1,\dots,u_n),v=(v_1,\dots,v_n)\in\R_+^n$, define the 
two-diagonal matrix 
\Eq{Auv}{
A(u,v):=
\left(\begin{array}{ccccc}
     0 & u_1  & \dots & 0 & 0\\
     v_1 & 0  & \dots & 0 & 0\\
     \vdots & \vdots & \ddots & \vdots & \vdots\\
     0 & 0 &  \dots & 0 & u_n  \\
     0 & 0 & \dots & v_n & 0
 \end{array}\right)\in\R_{+}^{(n+1)\times (n+1)}.
}
Then all the eigenvalues of $A(u,v)$ are real numbers. Furthermore, the eigenvalues of $A(u,v)$ are 
smaller than $1$ if and only if $w_1,\dots,w_n>0$, where $w_{-1}:=w_0:=1$, and 
\Eq{mu}{
  w_{k}:=w_{k-1}-u_{k}v_{k}w_{k-2} \qquad(k\in\{1,\dots,n\}).
}}

\begin{proof} In the sequel, denote by $I_k$ the unit matrix of the matrix algebra $\R^{k\times k}$, and for 
a square matrix $S\in\R^{k\times k}$, denote by $P_S$ the characteristic polynomial of $S$ defined for 
$\lambda\in C$ by $P_S(\lambda):=\det(\lambda I_k-S)$. 

Let $u=(u_1,\dots,u_n),v=(v_1,\dots,v_n)\in\R_+^n$ be fixed. Define $A_0(u,v):=0$ and 
\Eq{Akuv}{
A_k(u,v):=
\left(\begin{array}{ccccc}
     0 & u_1  & \dots & 0 & 0\\
     v_1 & 0  & \dots & 0 & 0\\
     \vdots & \vdots & \ddots & \vdots & \vdots\\
     0 & 0 &  \dots & 0 & u_k  \\
     0 & 0 & \dots & v_k & 0
 \end{array}\right)\in\R_{+}^{(k+1)\times (k+1)} \qquad(k\in\{1,\dots,n\}).
}
Then $A_n(u,v)=A(u,v)$. Observe that $P_{A_0(u,v)}(\lambda)=\lambda$ and 
$P_{A_1(u,v)}(\lambda)=\lambda^2-u_1v_1$. Expanding the determinant of the characteristic polynomial by its 
last row, we can easily deduce the following recursive formula for $P_{A_{k+1}(u,v)}$:
\Eq{rec}{
 P_{A_{k+1}(u,v)}(\lambda)
  =\lambda P_{A_k(u,v)}(\lambda)-u_{k+1}v_{k+1}P_{A_{k-1}(u,v)}(\lambda) \qquad(k\in\{1,\dots,n-1\}).
}

Now, we are going to prove that, for all $k\in\{1,\dots,n\}$, the characteristic polynomials of $A_k(u,v)$ 
and 
$A_k(\sqrt{uv},\sqrt{uv})$ are identical, where $\sqrt{uv}:=(\sqrt{u_1v_1},\dots,\sqrt{u_nv_n})$. 

We prove this statement by induction on $k$. For $k=0$, the statement is trivial. For $k=1$, we have that
\Eq{*}{
  P_{A_1(u,v)}(\lambda)=\lambda^2-u_1v_1=\lambda^2-\sqrt{u_1v_1}\sqrt{u_1v_1}
    =P_{A_1(\sqrt{uv},\sqrt{uv})}(\lambda).
}
Assume that we have established the identity $P_{A_j(u,v)}=P_{A_j(\sqrt{uv},\sqrt{uv})}$ for 
$j\leq k$. By using the recursive formula \eq{rec} twice and the inductive assumption, for 
$k\in\{1,\dots,n-1\}$, we get
\Eq{*}{
 P_{A_{k+1}(u,v)}(\lambda)
  &=\lambda P_{A_k(\sqrt{uv},\sqrt{uv})}(\lambda)
   -\sqrt{u_{k+1}v_{k+1}}\sqrt{u_{k+1}v_{k+1}}P_{A_{k-1}(\sqrt{uv},\sqrt{uv})}(\lambda)\\
  &=P_{A_{k+1}(\sqrt{uv},\sqrt{uv})}(\lambda).
}
This completes the proof of the identities $P_{A_k(u,v)}=P_{A_k(\sqrt{uv},\sqrt{uv})}$.

The matrix $A_n(\sqrt{uv},\sqrt{uv})$ is symmetric with real entries, therefore its characteristic polynomial 
has only real roots, whence it follows that the eigenvalues of $A_n(u,v)=A(u,v)$ are also real. 
The eigenvalues of $A_n(\sqrt{uv},\sqrt{uv})$ are smaller than one if and only if the eigenvalues of 
the symmetric matrix $I_{n+1}-A_n(\sqrt{uv},\sqrt{uv})$ are positive, which is equivalent to the positive 
definiteness of $I_{n+1}-A_n(\sqrt{uv},\sqrt{uv})$. In view of the Sylvester test, this holds if and only if 
all the leading principal minor determinants of $I_{n+1}-A_n(\sqrt{uv},\sqrt{uv})$ are positive, i.e., if
\Eq{pos}{
 P_{A_k(u,v)}(1)=P_{A_k(\sqrt{uv},\sqrt{uv})}(1)>0 \qquad(k\in\{0,\dots,n\}).
} 
By the recursive formula \eq{rec} applied for $\lambda=1$, it results that 
$P_{A_k(\sqrt{uv},\sqrt{uv})}(1)=w_k$ for all $k\in\{0,\dots,n\}$, therefore, \eq{pos} is equivalent to the
inequalities $w_1,\dots,w_n>0$. 
\end{proof}

The next result offers a sufficient condition in order that the eigenvalues of the matrix $A(u,v)$ be smaller 
than 1.

\Lem{suff}{
Let $n\in\N$ and $u=(u_{1},\dots, u_{n}), v=(v_{1},\dots, v_{n})\in\R^{n}$ be with positive components.
Assume that 
\Eq{uv}{
  v_1\leq 1, \qquad \max\{u_1+v_2,\dots,u_{n-1}+v_n\}\leq 1, \qquad u_n<1.
}
Then the system of inequalities $w_1,\dots,w_n>0$ holds, where $w_1,\dots,w_n$ are defined as in \thm{pd}.
Consequently, all the eigenvalues of the two-diagonal matrix $A(u,v)$ defined by \eq{Auv} are smaller than 1.}

\begin{proof} Observe that the positivity of $v_2,\dots,v_n$ and \eq{uv} yield that $u_1,\dots,u_n<1$.

To show that $w_k$ is positive for all $k\in\{1,\dots,n\}$, we shall prove that
\Eq{egy}{
   w_k>0 \qquad\mbox{and}\qquad (1-u_k)w_{k-1}\leq w_k<w_{k-1} \qquad(k\in\{1,\dots,n-1\}).
}
For $k=1$, the second chain of inequalities is equivalent to $1-u_1\leq 1-u_1v_1\leq 1$, which easily 
follows from $0<v_1\leq 1$ and $0<u_1$. Hence $w_1>0$ also holds. 

Assume that we have proved \eq{egy} for some $k\in\{1,\dots,n-1\}$. 
Then, using the recursion \eq{mu} and using the right hand side inequality in \eq{egy}, we get
\Eq{*}{
  w_{k+1}=w_k-u_{k+1}v_{k+1}w_{k-1}<w_k-u_{k+1}v_{k+1}w_k=w_k(1-u_{k+1}v_{k+1})<w_k.
}
On the other hand, using the upper estimate for $w_{k-1}$ obtained from \eq{egy}, we get
\Eq{*}{
  w_{k+1}=w_k-u_{k+1}v_{k+1}w_{k-1}
   &\geq w_k-u_{k+1}v_{k+1}\frac{w_k}{1-u_k}=w_k\frac{1-u_k-u_{k+1}v_{k+1}}{1-u_k}\\
   &\geq w_k\frac{1-u_k-u_{k+1}(1-u_k)}{1-u_k}=w_k(1-u_{k+1})>0,
}
which completes the proof of \eq{egy}.
\end{proof}

\Lem{2}{For $n\in\N$ and for the vectors $u=(u_1,\dots,u_n),v=(v_1,\dots,v_n)\in\R^n$ with positive 
components, define the two-diagonal matrix $A(u,v)$ by \eq{Auv}. Then there exists an eigenvector of $A(u,v)$ 
with positive components whose eigenvalue is also positive.}

\begin{proof} We follow the argument of the standard proof of the Perron--Frobenius Theorem. Consider the set 
\Eq{*}{
  S_{n+1}:=\{(x_0,\dots,x_{n})\in\R^{n+1}:x_0,\dots,x_n\geq0,\,x_0+\dots+x_n=1\}.
} 
Then $S_{n+1}$ is a compact convex set in $\R^{n+1}$. Let $u,v\in\R^n$ be fixed vectors with 
positive components and let $A_0,\dots,A_n$ be the row vectors of the matrix $A(u,v)$. Observe that 
\Eq{Au}{
  A(u,v)x=\big(\langle A_0,x\rangle,\dots,\langle A_n,x\rangle\big) \qquad(x\in\R^{n+1}),
} 
furthermore the sum $\langle A_0,x\rangle+\cdots+\langle A_n,x\rangle$ does not vanish on $S_{n+1}$. 
Indeed, if for some $x\in S_{n+1}$ we have $\langle A_0,x\rangle+\cdots+\langle A_n,x\rangle=0$, then, by 
the nonnegativity of the terms on the left hand side of this equation, it follows that 
$\langle A_{i},x\rangle=0$ for all $i\in\{0,\dots,n\}$. Using the positivity of the parameters $u_i$ and 
$v_i$, these equalities imply $x=0$, which contradicts $x\in S_{n+1}$. Now consider the mapping 
$F:S_{n+1}\to\R^{n+1}$ defined by
\Eq{*}{
  F(x):=\frac{A(u,v)x}{\langle A_0,x\rangle+\cdots+\langle A_n,x\rangle}\qquad(x\in S_{n+1}).
}
By \eq{Au}, we have that $F(S_{n+1})\subseteq S_{n+1}$, and $F$ is trivially continuous on $S_{n+1}$, hence,
by the Brouwer Fixed Point Theorem, there exists a fixed point $c\in S_{n+1}$ of the function $F$. Then we 
have
\Eq{*}{
  A(u,v)c=(\langle A_0,c\rangle+\dots+\langle A_n,c\rangle)F(c)
    =(\langle A_0,c\rangle+\dots+\langle A_n,c\rangle)c,
}
which shows that $c$ is an eigenvector of $A(u,v)$ with eigenvalue $\lambda:=\langle 
A_0,c\rangle+\cdots+\langle A_n,c\rangle>0$. Therefore, by $A(u,v)c=\lambda c$, the following system of 
equations hold for the coordinates $(c_0,\dots,c_n)$: 
\Eq{soe}{
\begin{array}{lclclcl}
u_1c_1 & & & = & \lambda c_0, \\
u_{i+1}c_{i+1} & + & v_{i}c_{i-1} & = & \lambda c_{i} \qquad (i\in\{1,\dots, n-1\}), \\
& & v_nc_{n-1} & = & \lambda c_n.
\end{array}
}
If $c_i=0$ for some $i\in\{0,\dots, n\}$, then the nonnegativity of the terms on the left hand side of the 
$i$th equation yields that $c_j=0$ for $j\in\{i-1,i+1\}\cap\{1,\dots, n-1\}$. This results that $c$ has to be 
zero, which contradicts $c\in S_{n+1}$.
\end{proof}

\section{Auxiliary results from fixed point theory}

For our purposes, we recall some notions and results related to fixed point theorems.

\Defn{*}{
We say that the function $d:X\times X\to\R$ is a \textit{semimetric} on the set $X$ if, for all 
$x,y\in X$, it possesses the following properties:
\begin{enumerate}\itemsep=0.2cm
\item\label{sm1} $d$ is \textit{positive definite}, i.e., for all $x,y\in X$, $d(x,y)\geq 0$ and $d(x,y)=0$ 
if and only if $x=y$,
\item\label{sm2} and $d$ is \textit{symmetric}, i.e., for all $x,y\in X$, $d(x,y)=d(y,x)$. 
\end{enumerate}
The pair $(X, d)$ is called \textit{semimetric space}.

If $(X, d_{X})$ and $(Y, d_{Y})$ are semimetric spaces then a function $f:X\to Y$ is called 
\textit{$L$-Lipschitzian} with respect to the pair of semimetrics $(d_{X}, d_{Y})$ if there exists $0\leq 
L$ such that
\Eq{L}{
d_{Y}(f(x), f(y))\leq Ld_{X}(x,y)\qquad(x,y\in X).
}
The function $f$ is said to have the \textit{Lipschitz property} if there exists $L\geq0$ such that \eq{L} 
holds. The \textit{Lipschitz modulus} of $f$ is defined by
\Eq{*}{
\lip f:=\sup\limits_{\substack{x,\,y\,\in\,X \\ x\,\neq\,y}}
        \frac{d_{Y}(f(x), f(y))}{d_{X}(x,y)}.
}
Obviously, $f$ has the Lipschitz property if and only if $\lip f$ is finite. If $\lambda:=\lip f<1$ then $f$ 
is called a \textit{$\lambda$-contraction}.}

It is an immediate consequence of these definitions, that for a subset $D\subseteq X$ and a contraction 
$f:D\to X$ with respect to the semimetric $d_{X}$, the map $f$ can have at most one fixed point in $D$. 
Indeed, if $x$ and $y$ are both fixed points of $f$  in $D$, then
\Eq{*}{
d_{X}(x,y)=d_{X}(f(x),f(y))\leq \lambda d_{X}(x,y),
}
which implies $d_{X}(x,y)\leq 0$, whence $x=y$ follows.

The following lemma is useful to compute the Lipschitz modulus of real valued functions.

\Lem{Lip}{Let $f,g:I\to\R$ be differentiable functions such that the derivative of $g$
does not vanish on the interval $I$. Then, for the Lipschitz modulus of the function $f\circ 
g^{-1}:g(I)\to\R$, the following formula holds:
\Eq{*}{
  \lip f\circ g^{-1}=\sup\limits_{t\,\in\,I}\bigg|\frac{f'(t)}{g'(t)}\bigg|.
}}

\begin{proof} Due to the assumptions of the lemma, $g:I\to g(I)$ is a continuous and strictly monotone 
function.
Therefore $g^{-1}:g(I)\to\R$ is well-defined. Thus, applying the Cauchy Mean Value Theorem, we have that
\Eq{*}{
  \lip f\circ g^{-1}
  = \sup\limits_{\substack{x,\,y\,\in\,g(I) \\ x\,\neq\,y}}
        \frac{|f\circ g^{-1}(x)-f\circ g^{-1}(y)|}{|x-y|}
  = \sup\limits_{\substack{u,\,v\,\in\,I \\ u\,\neq\,v}}
        \frac{|f(u)-f(v)|}{|g(u)-g(v)|}
  =\sup\limits_{t\,\in\,I}\bigg|\frac{f'(t)}{g'(t)}\bigg|.
}
\end{proof}

In what follows, we recall first the following generalization of the Tychonov Fixed Point Theorem established 
by Halpern and Bergman. For the formulation of this result, we define the notion of the \emph{inward set of a 
convex subset $K$} of a locally convex space $X$ by
\Eq{*}{
  \inw_K(x):=x+\R_+(K-x)=\{x+t(y-x)\mid y\in K,\,t\geq0\} \qquad(x\in K).
}
Observe that $K\subseteq\inw_K(x)$ holds for all $x\in K$. On the other hand, for an interior point $x$ of 
$K$, we have $\inw_K(x)=X$, therefore the inclusion $y\in \inw_K(x)$ is always trivial if $x\in 
K\setminus\partial K$, where the notation $\partial K$ stands for the set of boundary points of $K$.

\Thm{HB}{\textup{\textbf{(Halpern--Bergman)}}
Let $X$ be Hausdorff locally convex space and let $K\subseteq X$ be a compact convex set.
Let $f:K\to X$ be a continuous weakly inward map, i.e., assume that $f(x)\in \overline\inw_K(x)$ holds for 
all $x\in \partial K$. Then the set of fixed points of $f$ is a nonempty compact subset of $K$.}

If $f(K)\subseteq K$, then $f(x)\in \overline\inw_K(x)$ trivially holds for all $x\in \partial K$, therefore, 
in this case, the above result reduces to the Tychonov Fixed Point Theorem. 

The fixed point theorem stated below that we are going to use for the existence proofs in our main results is 
consequence of the Halpern--Bergman Fixed Point Theorem. It establishes the existence of the fixed 
point for continuous maps defined over a convex polyhedron.

\Thm{FPT}{
Let $c_1,\dots,c_m\in\R^n$ and $\gamma_1,\dots,\gamma_m\in\R$ and assume that the polyhedron
$K\subseteq\R^n$ defined by
\Eq{K}{
  K:=\big\{x\in\R^n\mid \langle c_k,x\rangle\leq\gamma_k,\,k\in\{1,\dots,m\}\big\}
}
is bounded. Let $f:K\to\R^n$ be a continuous function with the following property
\Eq{IC}{
   \langle c_k,f(x)\rangle\leq\gamma_k \qquad \mbox{for all $x\in K$ and for all $k\in\{1,\dots,m\}$
   such that $\langle c_k,x\rangle=\gamma_k$}.
}
Then the set of fixed points of $f$ is a nonempty compact subset of $K$.}

\begin{proof} By our assumption, $K$ is a compact convex set. It suffices to show that, for all $x\in K$, 
\Eq{II}{
  \inw_K(x)=\big\{u\in\R^n\mid \langle c_k,u\rangle\leq\gamma_k  
  \mbox{ for all $k\in\{1,\dots,m\}$ such that $\langle c_k,x\rangle=\gamma_k$}\big\}
}
because condition \eq{IC} then implies that $f(x)\in\inw_K(x)$ for all $x\in K$, whence the 
Halpern--Bergman Fixed Point Theorem yields the existence the fixed point of $f$.

Let $x\in K$ be fixed. If $u\in\inw_K(x)$, then there exists $y\in K$ and $t\geq0$ such that $u=(1-t)x+ty$.
Then, for $k\in\{1,\dots,m\}$ such that $\langle c_k,x\rangle=\gamma_k$, we have
\Eq{*}{
  \langle c_k,u\rangle = \langle c_k,(1-t)x+ty\rangle 
    = (1-t)\langle c_k,x\rangle + t\langle c_k,y\rangle
    = (1-t)\gamma_k + t\langle c_k,y\rangle \leq (1-t)\gamma_k + t\gamma_k=\gamma_k,
}
which proves the inclusion $\subseteq$ in \eq{II}.

For the reversed inclusion, let $u\in\R^n$ be an element such that $\langle c_k,u\rangle\leq\gamma_k$ 
for all $k\in\{1,\dots,m\}$ such that $\langle c_k,x\rangle=\gamma_k$. Choose $t>0$ such that
\Eq{*}{
  t\geq\frac{\langle c_k,u-x\rangle}{\gamma_k-\langle c_k,x\rangle} \qquad
  \mbox{for all $k\in\{1,\dots,m\}$ such that $\langle c_k,x\rangle<\gamma_k$}
}
and define $y\in\R^n$ by $y:=\tfrac{1}{t}(u-x)+x$. Then, distinguishing the cases whether $\langle 
c_k,x\rangle=\gamma_k$ or not, for every $k\in\{1,\dots,m\}$, we get that
\Eq{*}{
  \langle c_k,u-x\rangle \leq t (\gamma_k-\langle c_k,x\rangle).
}
Therefore, for every $k\in\{1,\dots,m\}$,
\Eq{*}{
  \langle c_k,y\rangle 
  = \langle c_k,\tfrac{1}{t}(u-x)+x\rangle
  \leq (\gamma_k-\langle c_k,x\rangle) + \langle c_k,x\rangle
  = \gamma_k.
}
This proves that $y\in K$. On the other hand, from the definition of $y$, we have that $u=(1-t)x+ty$, 
consequently, $u\in\inw_K(x)$.
\end{proof}

\section{The descendants of means}

Assume that we are given an $n\geq2$ member sequence of means $M_1,\dots,M_n:I_\leq^2\to\R$. 
In this section, we are going to deal with existence and uniqueness of an increasing sequence of means 
$N_1,\dots,N_n:I_\leq^2\to\R$ such that, for all $(x,y)\in I_\leq^2$, the identities 
\Eq{MN}{
N_1(x,y)&=M_1(x,N_2(x,y)),\\
N_i(x,y)&=M_i(N_{i-1}(x,y),N_{i+1}(x,y)) \qquad(i\in\{2,\dots,n-1\}),\\
N_n(x,y)&=M_n(N_{n-1}(x,y),y)
}
hold. The $i^\textrm{th}$ element of the sequence $N_1,\dots,N_n$ will be called the $i^\textrm{th}$
descendant of the $n$-tuple $(M_1,\dots,M_n)$ of means. Observe that \eq{MN} states that, for $(x,y)\in 
I_\leq^2$, the vector $(N_1(x,y),\dots,N_n(x,y))\in[x,y]_\leq^n$ is a fixed point of the mapping 
$\varphi_{(x,y)}:[x,y]_\leq^n\to\R^n$ defined by
\Eq{phid}{
\varphi_{(x,y)}(t_{1},\dots, t_n):=
\big(M_{1}(x,t_{2}), \dots, M_{i}(t_{i-1},t_{i+1}),\dots, M_n(t_{n-1},y)\big).
}

The first main result of this section establishes the existence and uniqueness of the fixed points of 
$\varphi_{(x,y)}$, i.e., the nonemptiness and singletonness of the set 
\Eq{Phid}{
\Phi_{(x,y)}:=
\big\{\xi\in[x,y]_\leq^n\,\big{|}\,\varphi_{(x,y)}(\xi)=\xi\big\}.
}
The existence and uniqueness of the fixed point is obvious if $x=y$, therefore, we restrict our attention to 
the case $x<y$.

\Thm{5p}{
Let $n\geq 2$ and $M_{1},\dots, M_n:I_\leq^2\to I$ be means.
For $(x,y)\in I_<^2$, define the mapping $\varphi_{(x,y)}$ and the fixed point set $\Phi_{(x,y)}$ 
by \eq{phid} and \eq{Phid}, respectively.  
Then, for all $(x,y)\in I_<^2$, the following statements hold:
\begin{enumerate}\itemsep=2mm
\item\label{5p1} If all the means $M_1,\dots, M_n$ are continuous, then the fixed point set 
  $\Phi_{(x,y)}$ is nonempty and compact. If the means $M_1,\dots, M_n$ are strict, 
	then $\Phi_{(x,y)}\subseteq\,]x, y[\,_<^n$.
\item\label{5p2} The set $\Phi_{(x,y)}$ is a singleton if there exist semimetrics 
$d_1,\dots,d_n:[x,y]^2\to\R_+$ 
such that the estimates
\Eq{lip}{
d_1(M_{1}(x,s),M_{1}(x,v))&\leq b_{1}d_2(s,v), \\
d_i(M_{i}(t,s),M_{i}(u,v))&\leq a_{i}d_{i-1}(t,u)+b_{i}d_{i+1}(s,v)\qquad (i\in\{2,\dots,n-1\}), \\
d_n(M_n(t,y),M_n(u,y))&\leq a_nd_{n-1}(t,u)
}
hold for all $t,s,u,v\in[x,y]$ with some positive real numbers $a_{2},\dots, a_n$ and $b_{1},\dots, 
b_{n-1}$ such that $w_1,\dots,w_{n-1}>0$, where $w_{-1}:=w_0:=1$ and 
\Eq{ab}{
  w_{i}:=w_{i-1}-a_{i+1}b_{i}w_{i-2}\qquad(i\in\{1,\dots,n-1\}).
}
\end{enumerate}}

\begin{proof}
Let $(x,y)\in I_<^2$ be arbitrarily fixed. Then the set $K:=[x,y]_\leq^n$ is a compact convex set which is 
characterized by the following $(n+1)$ inequalities: $(t_1,\dots,t_n)\in K$ holds if and only if
\Eq{ttt}{
  -t_1\leq -x,\qquad t_1-t_2\leq 0,\qquad\dots,\qquad t_{n-1}-t_n\leq 0,\qquad t_n\leq y.
}
Therefore, $K$ is a polyhedron of the form \eq{K} with $m=n+1$, suitably chosen 
vectors $c_1,\dots,c_{n+1}\in\R^n$ and scalars $\gamma_1,\dots,\gamma_{n+1}\in\R$. Thus, in order to show 
that the fixed point set of the continuous function $f:=\varphi_{x,y}$ is a nonempty compact subset of 
$K=[x,y]_\leq^n$, we need to verify that condition \eq{IC} is satisfied. 

For the sake of brevity, denote $t_0:=x$ and $t_{n+1}:=y$.
If, for some $k\in\{2,\dots,n\}$, the $k$th inequality holds with equality in \eq{ttt}, then $t_{k-1}=t_k$. 
Therefore, by the mean value property of the means $M_{k-1}$ and $M_k$, we get 
\Eq{*}{
  s_{k-1}=M_{k-1}(t_{k-2},t_k)\leq t_k=t_{k-1}\leq M_k(t_{k-1},t_{k+1})=s_k,
}
which proves that the vector $s$ satisfies the $k$th inequality in \eq{ttt}.

On the other hand, by the mean value properties of $M_1$ and $M_n$, we have $x\leq M_1(x,t_2)=s_1$ and 
$s_n=M_n(t_{n-1},y)\leq y$, therefore, $s$ also satisfies the first and last inequality in \eq{ttt} and thus 
the verification of condition \eq{IC} is complete.

To prove the second part of the statement (\ref{5p1}), assume that all the means $M_1,\dots, M_n$ are strict 
and let 
$(\xi_1,\dots,\xi_n)\in\Phi_{(x,y)}$. Then
\Eq{xi}{
  M_1(x,\xi_2)=\xi_1, \quad M_2(\xi_1,\xi_3)=\xi_2, \qquad \dots,\qquad M_n(\xi_{n-1},y)=\xi_n.
}
If $x=\xi_1$ then the strict mean property of $M_1$ and the identity $M_1(x,\xi_2)=\xi_1$ imply that 
$\xi_1=\xi_2$. Now the strict mean property of $M_2$ and the identity $M_2(\xi_1,\xi_3)=\xi_2$ yield that 
$\xi_2=\xi_3$.
Continuing this argument, it follows that $\xi_{n-1}=\xi_n$. Finally, the strict mean property of 
$M_n$ and $M_n(\xi_{n-1},y)=\xi_n$ imply that $\xi_n=y$. This leads to the contradiction 
$x=y$. Hence, we may assume that $x<\xi_1$. Applying the strict mean property of $M_1,\dots,M_n$ and the 
equalities in \eq{xi}, we get $\xi_{i}<\xi_{i+1}$ recursively for $i\in\{1,\dots,n-1\}$ and finally 
$\xi_n<y$, which proves that $(\xi_1,\dots, \xi_n)\in\,]x, y[\,_<^n$.

To prove (\ref{5p2}), assume that there exist semimetrics $d_1,\dots,d_n:[x,y]^2\to\R_+$ such that the 
estimates
in \eq{lip} hold and let $a:=(a_{2},\dots, a_n)$ and $b:=(b_{1},\dots, b_{n-1})$ such that each 
members of the sequence $w_1,\dots,w_{n-1}$, defined by \eq{ab} with $w_{-1}:=w_0:=1$, is positive. According 
to the 
previous lemmas, the matrix $A(a,b)$ has an eigenvector $c:=(c_{1},\dots, c_n)$ with positive components 
and with eigenvalue $0<\lambda<1$. This means that $c$ and $\lambda$ satisfy the following system of linear 
equations:
\Eq{abc}{
\begin{array}{lclclcl}
a_2c_2 & & & = & \lambda c_1, \\
a_{i+1}c_{i+1} & + & b_{i-1}c_{i-1} & = & \lambda c_{i} \qquad\quad (i\in\{2,\dots, n-1\}), \\
& & b_{n-1}c_{n-1} & = & \lambda c_n.
\end{array}
}

We show that $\varphi_{(x,y)}$ is a $\lambda$-contraction 
with respect to the semimetric $D_{c}:[x,y]^n\times[x,y]^n\to\R_{+}$ defined by
\Eq{*}{
D_{c}((u_{1},\dots, u_n),(v_{1},\dots, v_n))
:=c_{1}d_{1}(u_{1},v_{1})+\dots+c_nd_n(u_n,v_n)
}
for all $(u_{1},\dots, u_n),(v_{1},\dots, v_n)\in[x,y]^n$. To prove this, let $(t_{1},\dots, 
t_n)$ and $(s_{1},\dots, s_n)$ be arbitrary elements of $[x,y]_\leq^n$. For the sake of brevity, set 
$t_0=s_0=x$ and $t_n=s_n=y$. Using the estimates in \eq{lip} and then the identities in \eq{abc}, 
we obtain that
\Eq{*}{
D_{c}&(\varphi_{(x,y)}(t_{1},\dots, t_n),\varphi_{(x,y)}(s_{1},\dots, s_n))
 =\sum\limits_{i=1}^n c_{i}d_{i}\big(M_{i}(t_{i-1}, t_{i+1}),M_{i}(s_{i-1}, s_{i+1})\big) \\
&\leq c_{1}b_{1}d_{2}(t_{2},s_{2})
 +\bigg(\sum\limits_{i=2}^{n-1}c_{i}a_{i}d_{i-1}(t_{i-1},s_{i-1})+c_{i}b_{i+1}d_{i}(t_{i+1},s_{i+1})\bigg)
 + c_na_nd_n(t_{n-1},s_{n-1}) \\
&=\lambda\big(c_{1}d_{1}(t_{1},s_{1})+\dots+c_nd_n(t_n,s_n)\big)
 =\lambda D_{c}((t_{1},\dots, t_n),(s_{1},\dots, s_n)).
}
This results the uniqueness of the fixed point of $\varphi_{(x,y)}$.
\end{proof}

\Defn{SF}{
Let $n\geq 2$ and $M_{1},\dots, M_n:I_\leq^2\to\R$ be continuous means. For $i\in\{1,\dots,n\}$, we say that 
$N:I_\leq^2\to\R$ is an \textit{$i^{\mbox{\tiny\rm th}}$ 
descendant of the $n$-tuple of means $(M_{1},\dots,M_n)$} if, for all $(x,y)\in I_\leq^2$, we have
\Eq{sf}{
N(x,y)\in\bigcup\big\{\xi_i\mid(\xi_1,\dots,\xi_n)\,\in\,\Phi_{(x,y)}\big\}
\qquad\text{if }x<y \qquad \text{ and }\qquad 
N(x,y)=x\quad\text{if }x=y,
}
where $\Phi_{(x,y)}$ is the fixed point set of the mapping $\varphi_{(x,y)}:[x,y]_\leq^n\to\R^n$ defined by 
\eq{phid}. The class of all such functions is denoted by $\D_i(M_{1},\dots, M_n)$.}

Note that, in view of \thm{5p}, the continuity of the means $M_{1},\dots, M_{n}$ implies that the descendant 
functions are well-defined. As a direct consequence of the compactness of the fixed point 
set $\Phi_{(x,y)}$, we obtain that the family $\D_i(M_{1},\dots, M_{n})$ has a minimal and a maximal 
element in the following sense: there exist $N^-_i,N^+_i\in\D_i(M_{1},\dots, M_{n})$  
such that $N^-_i(x,y)\leq N(x,y)\leq N^+_i(x,y)$) for all $N\in\D_i(M_{1},\dots, M_{n})$ 
and for all $x,y\in I$. It is also obvious that each element of 
$\D_i(M_{1},\dots, M_n)$ is a strict mean provided that all the means $M_{1},\dots, M_n$ are strict.

\Rem{1}{The uniqueness of the fixed point of the map $\varphi_{(x,y)}$ cannot be stated in general. For 
instance, let $n\geq2$, and let $M_1:=\max$, $M_n:=\min$ and $M_i:=\AM_{\frac12}$ for $i\in\{2,\dots,n-1\}$ 
over the interval $\R$. Then, for $(x,y)\in\R_<^2$, the fixed point equation 
$(t_1,\dots,t_n)=\varphi_{(x,y)}(t_1,\dots,t_n)$ is equivalent to
\Eq{*}{
  (t_1,\dots,t_n)=\Big(t_2,\frac{t_1+t_3}{2},\dots,\frac{t_{n-2}+t_n}{2},t_{n-1}\Big).
}
An easy computation shows that this equality is satisfied if and only if $t_1=\dots=t_n$. Therefore,
$\Phi_{(x,y)}=\{(t_1,\dots,t_n)\mid t_1=\dots=t_n\in[x,y]\}$.}

Considering Matkowski means, we obtain useful corollaries of \thm{5p}.

\Thm{tmat}{
Let $n\geq 2$ and $f_1,\dots,f_n,g_1,\dots,g_n:I\to\R$ be continuous, 
strictly increasing functions. For $(x,y)\in I_<^2$, define the function 
$\varphi_{(x,y)}:[x,y]_\leq^n\to\R^n$ as
\Eq{tmat}{
\varphi_{(x,y)}(t_{1},\dots, t_n):=
\big(\MM_{f_1,\,g_1}(x,t_{2}), \dots, \MM_{f_i,\,g_i}(t_{i-1},t_{i+1}),\dots, \MM_{f_n,\,g_n}(t_{n-1},y)\big).
}
Then, for $(x,y)\in I_<^2$, the fixed point set $\Phi_{(x,y)}$ defined by \eq{Phid}
is nonempty and compact. Furthermore, $\Phi_{(x,y)}$ is a singleton if
\Eq{lipc}{
 a_i&:=\lip\big[f_{i}\circ(f_{i-1}+g_{i-1})^{-1}\big]<+\infty \qquad(i\in\{2,\dots,n\}),\\[1mm]
 b_i&:=\lip\big[g_{i}\circ(f_{i+1}+g_{i+1})^{-1}\big]<+\infty \qquad(i\in\{1,\dots,n-1\}),
}
and if the constants $w_1,\dots,w_{n-1}$ defined by \eq{ab} are positive.
}

\begin{proof} Because of the definition of the means $\MM_{f_{1},\,g_{1}},\dots, \MM_{f_n,\,g_n}$, for all 
$(x,y)\in I_<^2$, the mapping $\varphi_{(x,y)}$ is continuous, thus, based on the \thm{5p}, the 
corresponding fixed point set $\Phi_{(x,y)}$ is a nonempty compact subset of $[x,y]_\leq^n$. Due to the 
strictness of generalized quasi-arithmetic means it also follows that $\Phi_{(x,y)}\subseteq\,]x, y[\,_<^n$.

Now assume that \eq{lipc} and $w_1,\dots,w_{n-1}>0$ hold and fix a point $(x,y)\in I_<^2$. To show that 
$\Phi_{(x,y)}$ is a singleton, for $i\in\{1,\dots, n\}$, define the semimetrics $d_i:I\times I\to\R_{+}$ 
as
\Eq{*}{
d_i(s,t):=
|(f_i+g_i)(s)-(f_i+g_i)(t)|\qquad(s,t\in I).
}
Note that in our case, for all $i\in\{1,\dots, n\}$, the function $d_i$ is a \textit{metric}, i.e., in 
addition 
of the properties (\ref{sm1}) and (\ref{sm2}) of semimetrics, $d_i$ also satisfies the \textit{triangle 
inequality}, namely, for all $i\in\{1,\dots, n\}$, we have
\Eq{*}{
d_i(s,t)\leq d_i(s,r)+d_i(r,t)\qquad(r,s,t\in I).
}
Let $t,s,u,v\in [x,y]$ be arbitrary. Then, for all $i\in\{2,\dots,n-1\}$, we have the following estimation:
\Eq{*}{
d_i(M_i(t,s), M_i(u,v))
&=|(f_i+g_i)(\MM_{f_i,\,g_i}(t,s))-(f_i+g_i)(\MM_{f_i,\,g_i}(u,v))|\\
&=|f_i(t)+g_i(s)-f_i(u)-g_i(v)|\leq|f_i(t)-f_i(u)|+|g_i(s)-g_i(v)|\\
&\leq\lip\big[f_i\circ(f_{i-1}+g_{i-1})^{-1}\big]d_{i-1}(t,u)
  +\lip\big[g_i\circ(f_{i+1}+g_{i+1})^{-1}\big]d_{i+1}(s,v) \\
&=a_id_{i-1}(t,u)+b_id_{i+1}(s,v).
}
On the other hand, for $i=1$ and $i=n$, we get
\Eq{*}{
d_{1}(M_{1}(x,s), M_{1}(x,v))&\leq b_{1}d_{2}(s,v)
\qquad\text{and}\qquad
d_n(M_n(t,y), M_n(u,y))&\leq a_nd_{n-1}(t,u).
}
Therefore, all the estimates in \eq{lip} are satisfied. Thus, in view of the \thm{5p}, for all 
$(x,y)\in I_<^2$, the fixed point set $\Phi_{(x,y)}$ is indeed a singleton.
\end{proof}

\Cor{ccmat}{
Let $n\geq 2$ and $f_1,\dots,f_n,g_1,\dots,g_n:I\to\R$ be differentiable, strictly increasing 
functions such that $(f_i+g_i)'$ does not vanish on $I$ for all $i\in\{1,\dots, n\}$. Assume further
that
\Eq{lipab}{
a_i:=\sup_{t\,\in\,I}\big[f'_i\cdot(f'_{i-1}+g'_{i-1})^{-1}\big](t)&<+\infty \qquad(i\in\{2,\dots,n\}),\\
b_i:=\sup_{t\,\in\,I}\big[g'_i\cdot(f'_{i+1}+g'_{i+1})^{-1}\big](t)&<+\infty \qquad(i\in\{1,\dots,n-1\}).
}
Finally, for $(x,y)\in I_<^2$, define the function $\varphi_{(x,y)}:[x,y]_\leq^n\to\R^n$ as in \eq{tmat}. 
Then, for all $(x,y)\in I_<^2$, the fixed point set $\Phi_{(x,y)}$ defined by \eq{Phid} is a nonempty 
compact subset of 
$[x,y]_\leq^n$, and, it is a singleton if the constants $w_1,\dots,w_{n-1}$ defined by \eq{ab} are positive.}

\begin{proof}
In view of \thm{tmat}, we only need to verify that $\Phi_{(x,y)}$ is a singleton, which in turn is obvious. 
Using \lem{Lip} and the conditions in \eq{lipab}, one can easily see that the estimations in \eq{lipc} of 
\thm{tmat} hold, i.e., the constants $a_{2},\dots, a_{n}$ and $b_{1},\dots, b_{n-1}$
are real numbers. Thus the proof is complete.
\end{proof}

\Thm{cmat}{
Let $n\geq 2$, $s_{1},\dots, s_n\in\,]0,1[\,$, and $h:I\to\R$ be a continuous, strictly increasing function.
Then, for all all $(x,y)\in I_<^2$, the fixed point set $\Phi_{(x,y)}$ (defined by \eq{Phid}) of the mapping 
$\varphi:[x,y]_\leq^n\to\R^n$ defined by
\Eq{*}{
\varphi_{(x,y)}(t_{1},\dots, t_n):=
\big(\MM_{s_{1}h,\,(1-s_{1})h}(x,t_{2}), \dots, \MM_{s_{i}h,\,(1-s_{i})h}(t_{i-1},t_{i+1}),\dots, 
\MM_{s_nh,\,(1-s_n)h}(t_{n-1},y)\big)
}
is the singleton 
$\big\{\big(\MM_{\sigma_{1}h,\,(1-\sigma_{1})h}(x,y),\dots,\MM_{\sigma_nh,\,(1-\sigma_n)h}(x,y)\big)\big\}$, 
where
\Eq{si}{
 \sigma_i:=\bigg(\sum_{j=i}^n\prod_{k=1}^{j}\frac{s_k}{1-s_k}\bigg)
   \bigg(\sum_{j=0}^n\prod_{k=1}^{j}\frac{s_k}{1-s_k}\bigg)^{-1} \qquad(i\in\{1,\dots,n\}).
}}

\begin{proof}
In order to apply \thm{tmat}, let $f_i:=s_i\cdot h$ and $g_i:=(1-s_i)\cdot h$ for 
$i\in\{1,\dots,n\}$. Then it immediately follows that the fixed point set 
$\Phi_{(x,y)}$ is nonempty and compact for all $(x,y)\in I_<^2$. 

To show that $\Phi_{(x,y)}$ is a singleton define the constants $a_{2},\dots, a_n$, 
$b_{1},\dots, b_{n-1}$, and $w_1,\dots,w_{n-1}$ as in \thm{tmat}. We need to show that conditions \eq{lipc} 
and 
$w_1,\dots,w_{n-1}>0$ hold. Observe that, for $i\in\{1,\dots, n\}$, we have $f_i+g_i=h$ and
\Eq{*}{
a_i&=\lip\big[f_i\circ(f_{i-1}+g_{i-1})^{-1}\big]
 =\lip [s_i\cdot h\circ h^{-1}]=s_i 
\qquad(i\in\{2,\dots,n\}), \\
b_i&=\lip\big[g_i\circ(f_{i+1}+g_{i+1})^{-1}\big]
 =\lip [(1-s_i)\cdot h\circ h^{-1}]=1-s_i 
\qquad(i\in\{1,\dots,n-1\}).
}
Thus each of the constants $a_{2},\dots, a_n$ and $b_{1},\dots, b_{n-1}$ are finite, on the other hand, 
under the notation $(u_1,\dots,u_{n-1}):=(a_2,\dots,a_n)$ and $(v_1,\dots,v_{n-1}):=(b_1,\dots,b_{n-1})$, 
they satisfy the condition \eq{uv} of \lem{suff}. Therefore, the inequalities $w_1,\dots,w_{n-1}>0$ hold
and hence $\Phi_{(x,y)}$ is a singleton.

Finally, we verify that, for all $(x,y)\in I_<^2$, the vector
$\big(\MM_{\sigma_{1}h,\,(1-\sigma_{1})h}(x,y),\dots,\MM_{\sigma_nh,\,(1-\sigma_n)h}(x,y)\big)$ is a fixed 
point of $\varphi_{(x,y)}$. For this purpose, we show first that $\sigma_1,\dots,\sigma_n$ fulfill the 
following system of linear equations:
\Eq{sig}{
  \sigma_1&=s_1+(1-s_1)\sigma_2,\\
  \sigma_i&=s_i\sigma_{i-1}+(1-s_i)\sigma_{i+1} \qquad(i\in\{2,\dots,n-1\}),\\
  \sigma_n&=s_n\sigma_{n-1}.
}
We prove the above equality for $i\in\{2,\dots,n-1\}$. First observe that
\Eq{*}{
  \prod_{k=1}^{i}\frac{s_k}{1-s_k}
  =\frac{s_i}{1-s_i}\prod_{k=1}^{i-1}\frac{s_k}{1-s_k}
  =s_i\Big(1+\frac{s_i}{1-s_i}\Big)\prod_{k=1}^{i-1}\frac{s_k}{1-s_k}
  =s_i\bigg(\prod_{k=1}^{i-1}\frac{s_k}{1-s_k}+\prod_{k=1}^{i}\frac{s_k}{1-s_k}\bigg).
}
Adding this identity to the equality 
\Eq{*}{
  \sum_{j=i+1}^n\prod_{k=1}^{j}\frac{s_k}{1-s_k}
   =s_i\sum_{j=i+1}^n\prod_{k=1}^{j}\frac{s_k}{1-s_k}
    +(1-s_i)\sum_{j=i+1}^n\prod_{k=1}^{j}\frac{s_k}{1-s_k} 
}
side by side, we get the desired identity $\sigma_i=s_i\sigma_{i-1}+(1-s_i)\sigma_{i+1}$.
In the cases $i=1$ and $i=n$ the proof of \eq{sig} is completely analogous.

Using \eq{sig}, after some calculation we easily get that
\Eq{*}{
\MM_{\sigma_{1}h,\,(1-\sigma_{1})h}(x,y)
 &=\MM_{s_{1}h,\,(1-s_{1})h}(x,\MM_{\sigma_{2}h,\,(1-\sigma_{2})h}(x,y)), \\
\MM_{\sigma_{i}h,\,(1-\sigma_{i})h}(x,y)
&=\MM_{s_{i}h,\,(1-s_{i})h}(\MM_{\sigma_{i-1}h,\,(1-\sigma_{i-1})h}(x,y),\MM_{\sigma_{i+1}h,\,(1-\sigma_{i+1}
)h} (x ,
y))   
 \quad(i\in\{2,\dots,n-1\}),\\
\MM_{\sigma_nh,\,(1-\sigma_n)h}(x,y)   
 &=\MM_{s_nh,\,(1-s_n)h}(\MM_{\sigma_{n-1}h,\,(1-\sigma_{n-1})h}(x,y),y),
}
which proves that 
$\big(\MM_{\sigma_{1}h,\,(1-\sigma_{1})h}(x,y),\dots,\MM_{\sigma_nh,\,(1-\sigma_n)h}(x,y)\big)$ is indeed a 
fixed point of $\varphi_{(x,y)}$.
\end{proof}

\Thm{rmat}{
Let $n\geq 2$, $j\in\{1,\dots, n\}$ and $p,q, h_{1},\dots, h_{n-1}:I\to\R$ be 
continuous, strictly increasing functions, and set $h_0:=h_n:=0$. For $(x,y)\in I_<^2$, define the 
mapping $\varphi_{(x,y)}:[x,y]_\leq^n\to\R^n$ by \eq{phid}, where
\Eq{*}{
M_{i}:=
\begin{cases}
\MM_{p+h_{i-1},\,h_i}, & \text{if } i\in\{1,\dots, j-1\}, \\[2mm]
\MM_{p+h_{i-1},\,h_{i}+q} & \text{if } i=j, \\[2mm]
\MM_{h_{i-1},\,h_{i}+q}, & \text{if } i\in\{j+1,\dots, n\}.
\end{cases}
}
Then, for $(x,y)\in I_<^2$, the fixed point set $\Phi_{(x,y)}$ defined by \eq{Phid}
is the singleton $\{(\xi_1,\dots,\xi_n)\}$, where the coordinates are defined by the following
two-way recursion:
\Eq{xij}{
\xi_j:=\MM_{p,\,q}(x,y) \qquad\mbox{and}\qquad
\xi_i:=
\begin{cases}
\MM_{p,\,h_{i}}(x,\xi_{i+1}) &\text{if }i\in\{1,\dots, j-1\}, \\[2mm]
\MM_{h_{i-1},\,q}(\xi_{i-1},y) &\text{if }i\in\{j+1,\dots, n\}.
\end{cases}
}}

\begin{proof} Let $(x,y)\in I_<^2$ be fixed. By \thm{5p}, the set $\Phi_{(x,y)}$ is nonempty.
Let $(\xi_1,\dots,\xi_n)\in \Phi_{(x,y)}$ be arbitrary and denote $\xi_0:=x$ and $\xi_{n+1}:=y$. Then, by the 
definition of Matkowski means, we have
\Eq{sys}{
  (p+h_{i-1}+h_i)(\xi_i)&=(p+h_{i-1})(\xi_{i-1})+h_i(\xi_{i+1}), && \qquad \text{if } i\in\{1,\dots, j-1\}, 
\\[1mm]
  (p+h_{i-1}+h_{i}+q)(\xi_i)&=(p+h_{i-1})(\xi_{i-1})+(h_{i}+q)(\xi_{i+1}), && \qquad \text{if } i=j, \\[1mm]
  (h_{i-1}+h_{i}+q)(\xi_i)&=h_{i-1}(\xi_{i-1})+(h_{i}+q)(\xi_{i+1}), && \qquad \text{if } i\in\{j+1,\dots,n\}.
}
Adding up these inequalities for $i\in\{1,\dots,n\}$ side by side, it follows that
\Eq{*}{
  p(\xi_j)+h_0(\xi_1)+h_n(\xi_n)+q(\xi_j)=p(\xi_0)+h_0(\xi_0)+h_n(\xi_{n+1})+q(\xi_{n+1}).
}
This simplifies to 
\Eq{*}{
  (p+q)(\xi_j)=p(x)+q(y),
}
which is equivalent to the equality on the left hand side of \eq{xij}. This computation also shows that 
$\xi_j$ is uniquely determined.

To prove the first equality on the right hand side of \eq{xij}, assume that $1\leq j-1$ and let 
$k\in\{1,\dots, j-1\}$ be fixed. Adding up the equalities in \eq{sys} for $i\in\{1,\dots,k\}$, we arrive at
\Eq{*}{
  p(\xi_k)+h_0(\xi_1)+h_k(\xi_k)=p(\xi_0)+h_0(\xi_0)+h_k(\xi_{k+1}),
}
which reduces to $(p+h_k)(\xi_k)=p(x)+h_k(\xi_{k+1})$ proving the first equality on the right hand side 
of \eq{xij} for $i=k$.

Analogously, to verify the second equality on the right hand side of \eq{xij}, assume that $j+1\leq n$ and
let $k\in\{j+1,\dots, n\}$ be fixed. Adding up the equalities in \eq{sys} for $i\in\{k,\dots,n\}$, we obtain
\Eq{*}{
  h_{k-1}(\xi_k)+h_n(\xi_n)+q(\xi_k)=h_{k-1}(\xi_{k-1})+h_n(\xi_{n+1})+q(\xi_{n+1}).
}
This yields $(h_{k-1}+q)(\xi_k)=h_{k-1}(\xi_{k-1})+q(y)$, which validates the second equality on the 
right hand side of \eq{xij} for $i=k$.

In view of the uniqueness of $\xi_j$ and the recursive system of equalities on the right hand side of 
\eq{xij}, we can see that, for $i\neq j$, the value of $\xi_i$ is also uniquely determined.
\end{proof}

\section{Upper- and lower second-order divided differences}

Consider the following binary operations on the extended real line $\RX:=\R\cup\{-\infty,+\infty\}$: for two 
extended real numbers $x,y$, their \textit{upper} and \textit{lower sums} are defined by
\Eq{*}{
x\uplus y:=
\begin{cases}
x+y, &\text{if } \max\{x,y\}<+\infty, \\
+\infty, & \text{if } \max\{x,y\}=+\infty,
\end{cases} 
\qquad
x\lplus y:=
\begin{cases}
x+y, &\text{if } \min\{x,y\}>-\infty, \\
-\infty, & \text{if } \min\{x,y\}=-\infty,
\end{cases}
}
respectively.
It is easy to see, that the pairs $(\RX,\uplus )$ and $(\RX,\lplus )$ are commutative semigroups. 
Apart from the standard cases, the only difference between these operations is that 
\Eq{*}{
  (-\infty)\uplus(+\infty)=(+\infty)\uplus(-\infty)=+\infty
   \qquad\mbox{and}\qquad
  (-\infty)\lplus(+\infty)=(+\infty)\lplus(-\infty)=-\infty.
}
Furthermore, the both of the operations $\uplus$ and $\lplus$ restricted to pairs of real numbers are 
the same as the standard addition of the reals. As direct
consequences of the definitions, for all $x,y\in\RX$, we have the following easy-to-see properties:
\Eq{pop}{
x\lplus y\leq x\uplus y \qquad\mbox{and}\qquad -(x\lplus y)=(-x)\uplus (-y),
}
furthermore, we have the following equivalences:
\Eq{pdd}{
  &0\leq x\uplus y \quad\Leftrightarrow\quad -x\leq y \qquad\mbox{and}\qquad
   0\leq x\lplus y \quad\Leftrightarrow\quad \big(-\infty<\min\{x,y\}\mbox{ and }-x\leq y\big), \\
  &x\lplus y \leq0 \quad\Leftrightarrow\quad x\leq -y \qquad\mbox{and}\qquad
   x\uplus y \leq0 \quad\Leftrightarrow\quad \big(\max\{x,y\}<+\infty\mbox{ and }x\leq -y\big). \\
}

\Defn{dd}{
Let $D\subseteq\R$ and $f:D\to\RX$. The \textit{upper second-order divided difference} of $f$ at three 
distinct points $x,y,z$ of $D$ is an extended real number defined by
\Eq{*}{
 \left\lceil x,y,z;f\right\rceil
 :=\frac{f(x)}{(y-x)(z-x)}\dot{+}\frac{f(y)}{(x-y)(z-y)}\dot{+}\frac{f(z)}{(x-z)(y-z)}.
}
Similarly, the \textit{lower second-order divided difference} of $f$ at the distinct points $x,y,z$
of $D$ is
\Eq{*}{
 \left\lfloor x,y,z;f\right\rfloor
 :=\frac{f(x)}{(y-x)(z-x)}\lplus \frac{f(y)}{(x-y)(z-y)}\lplus\frac{f(z)}{(x-z)(y-z)}.
}}
Obviously, the above second-order divided differences are symmetric functions of $(x,y,z)$.
Observe that if the inequalities $x<y<z$ hold, then the coefficients of $f(x)$ and $f(z)$ are positive and 
the coefficient to $f(z)$ is negative.

As a direct consequence of the above definition and \eq{pop} we obtain

\Prp{1p}{
Let $D\subseteq\R$ and $f:D\to\RX$. Then, for all distinct points $x<y<z$ of $D$,
\Eq{*}{
\lfloor x,y,z;f \rfloor\leq \lceil x,y,z;f \rceil
\qquad\text{and}\qquad
-\lfloor x,y,z;f\rfloor=\lceil x,y,z; -f\rceil.
}}

\Prp{eci}{
\textup{(Extended Chain Inequality)}
Let $D\subseteq\R$ and $f:D\to\RX$. Then, for all $n\in\N$ and $x_{0}<x_{1}<\dots<x_{n+1}$ 
in $D$ and for all $i\in\{1,\dots, n\}$ the following inequalities hold:
\Eq{*}{
\min\limits_{1\,\leq\,j\,\leq\,n}\lfloor x_{j-1},x_{j},x_{j+1};f\rfloor
\leq\lfloor x_{0},x_{i},x_{n+1};f\rfloor
\leq\lceil x_{0},x_{i},x_{n+1};f\rceil
\leq\max\limits_{1\,\leq\,j\,\leq\,n}\lceil x_{j-1},x_{j},x_{j+1};f\rceil.
}}

\begin{proof}
We only need to prove the first inequality, because the second one is trivial and the 
last one is the consequence of the first and \prp{1p}. 

The statement is trivial for $n=1$, therefore we may assume that $n\geq 2$. Let $x_{0}<x_{1}<\dots<x_{n+1}$ 
be 
arbitrary elements of $D$ and $i\in\{1,\dots, n\}$. If either the left hand side of the first inequality 
equals 
$-\infty$ or the right hand side equals 
$+\infty$, then there is nothing to prove. In the remaining case, for all $j\in\{1,\dots, n\}$, we have 
that $\lfloor x_{j-1},x_{j},x_{j+1};f\rfloor>-\infty$ and $\lfloor x_{0},x_{i},x_{n+1};f\rfloor<+\infty$.
The first inequality implies, for all $j\in\{1,\dots, n\}$ that
\Eq{*}{
\min\{f(x_{j-1}),-f(x_{j}),f(x_{j+1})\}>-\infty.
}
In view of $n\geq 2$, the set $\{1,\dots, n\}$ contains at least two elements, therefore, for all 
$j\in\{1,\dots,n\}$, we obtain that $f(x_j)\in\R$ and $\min\{f(x_0),f(x_{n+1})\}>-\infty$. Thus, $f(x_i)\in\R$
and hence the inequality $\lfloor x_{0},x_{i},x_{n+1};f\rfloor<+\infty$ yields 
$\max\{f(x_0),f(x_{n+1})\}<+\infty$,
which proves that, for all $j\in\{0,\dots,n+1\}$, we have $f(x_j)\in\R$. In this case, the first inequality 
is a consequence of 
\cite[Corollary 1]{NikPal03}.
\end{proof}

\section{Upper and lower $M$-convexity}
 
\Defn{mc}{
For a fixed strict mean $M:I_\leq^2\to\R$, we say that the function $f:I\to\RX$ is 
\textit{lower $M$-convex} if
\Eq{lmc}{
\left\lfloor x,M(x,y),y;f\right\rfloor\geq 0\qquad \big((x,y)\in I_<^2\big)
}
holds. On the other hand, the function $f$ is called \textit{upper $M$-convex} if
\Eq{umc}{
\left\lceil x,M(x,y),y;f\right\rceil\geq 0
}
holds on the same domain.}

Note that, due to the property \eq{pop} if $f$ is lower $M$-convex, then it is also 
upper $M$-convex.

The \textit{lower} and \textit{upper $M$-concavity} of functions can be also 
interpreted, namely we may consider \eq{lmc} and \eq{umc} with the reverse inequality. It is easy to check, 
that these definitions are equivalent to the upper and lower $M$-convexity of 
the function $-f$, respectively.

\Lem{EFC}{
Let $M:I_\leq^2\to\R$ be a strict mean and $f:I\to\RX$. Then the following statements hold.
\begin{enumerate}[(a)]\itemsep=0.2cm
\item The function $f$ is lower $M$-convex if and only if
$f(u)>-\infty$ for all $u\in I$ and, for all $(x,y)\in I_<^2$, the inequalities $f(M(x,y))<+\infty$ and 
\Eq{lmc+}{
f(M(x,y))\leq\frac{y-M(x,y)}{y-x}f(x)+\frac{M(x,y)-x}{y-x}f(y)
}
hold.
\item The function $f$ is upper $M$-convex if and only if, for all 
$(x,y)\in I_<^2$, the inequality
\Eq{umc+}{
f(M(x,y))\leq\frac{y-M(x,y)}{y-x}f(x)\uplus\frac{M(x,y)-x}{y-x}f(y)
}
holds.
\end{enumerate}
}

\begin{proof}
First we prove the statement (b). Suppose that $f$ is upper $M$-convex,
which means $\left\lceil x,M(x,y),y;f\right\rceil\geq 0$ for all $(x,y)\in I_<^2$. Due 
to the first property of upper addition in \eq{pdd}, this inequality is equivalent 
to
\Eq{6+}{
\frac{f(M(x,y))}{(M(x,y)-x)(y-M(x,y))}\leq
\frac{f(x)}{(M(x,y)-x)(y-x)}\uplus\frac{f(y)}{(x-y)(M(x,y)-y)},
}
where $(x,y)\in I_<^2$. Using that $(M(x,y)-x)(y-M(x,y))$ is positive, we obtain, for all 
$(x,y)\in I_<^2$, that \eq{umc+} is valid. 

To prove the reverse implication of (b), suppose that \eq{umc+} holds on the domain indicated. Then
\eq{6+} is also valid and, in view of the first property of upper addition in \eq{pdd}, this implies 
\eq{umc+}.

In the second step we prove the statement (a). Suppose that $f$ is 
lower $M$-convex, i.e. we have $\left\lfloor x,M(x,y),y;f\right\rfloor\geq 0$ for all $(x,y)\in I_<^2$.
Due to the first property of lower addition in \eq{pdd}, it follows that, $(x,y)\in I_<^2$, we have 
$-\infty<\min\{f(x), -f(M(x,y)), f(y)\}$ and
\Eq{6-}{
\frac{f(M(x,y))}{(M(x,y)-x)(y-M(x,y))}\leq
\frac{f(x)}{(M(x,y)-x)(y-x)}\lplus\frac{f(y)}{(x-y)(M(x,y)-y)},
}
Thus, for all $u\in I$, we get $-\infty<f(u)$ and, by the positivity of $(M(x,y)-x)(y-M(x,y))$, \eq{6-}
is equivalent to \eq{lmc+} and $f(M(x,y))<+\infty$ on the domain indicated.

To prove the reversed implication of the statement (a), suppose that $f(M(x,y))<+\infty$ and \eq{lmc+} hold 
for all $(x,y)\in I_<^2$ and we have $-\infty<f(u)$ for all $u\in I$. Then \eq{6-} is also valid and, 
in view of the first property of lower addition in \eq{pdd}, this implies \eq{lmc+}.
\end{proof}

In the following proposition we show that, for certain rational numbers $t$, there exists an upper 
$\AM_t$-convex extended real valued function $f$, which is not upper $\AM_{1-t}$-convex. Therefore, $f$ is 
not $t$-convex. It is an open problem if there exists a \textit{real-valued} function $f$ with these 
properties. This result is analogous to that of Lewicki and Olbry\'s \cite{LewOlb14} (which works for 
transcendental values of $t$).

\Prp{x1}{
Denote by $\Q_0$ and $\Q_1$ the following subsets of the rationals:
\Eq{*}{
  \Q_0:=\Big\{\frac{2k}{2n-1}\,\Big|\, k\in\Z,\,n\in\N\Big\} \qquad\mbox{and}\qquad
  \Q_1:=\Big\{\frac{2k-1}{2n-1}\,\Big|\, k\in\Z,\,n\in\N\Big\}.
}
Then $\Q_0$ and $\Q_1$ are disjoint subsets of $\Q$ and
\Eq{QQ}{
  \Q_0+\Q_0&\subseteq \Q_0,\qquad & \Q_0+\Q_1&\subseteq \Q_1,\qquad & \Q_1+\Q_1&\subseteq \Q_0, \\
  \Q_0\Q_0&\subseteq \Q_0,\qquad & \Q_0\Q_1&\subseteq \Q_0,\qquad & \Q_1\Q_1&\subseteq \Q_1.
}
Let $I\subseteq\R$ be an interval such that $a:=\sup I\in I\cap\Q_1$.
Let $h:I\to\R$ be an arbitrary convex function and define the function $f:I\to\RX$ by
\Eq{*}{
f(x):=
\begin{cases}
h(x) & \text{if }x\in (I\cap\Q_0)\cup\{a\}, \\[2mm]
+\infty & \text{if }x\in I\setminus (\Q_0\cup\{a\}).
\end{cases}
}
Then, for all $t\in\,]0,1[\,\cap\,\Q_1$, the function $f$ is upper $\AM_{t}$-convex and 
is not upper $\AM_{1-t}$-convex.
}

\begin{proof}
The inclusions in \eq{QQ} follow from elementary calculation with rational fractions.

Let $x,y\in I$ with $x<y$ and $t\in\,]0,1[\,\cap\,\Q_1$ be arbitrarily fixed. Then $1-t\in\Q_0$. We need to 
check that \eq{umc}
is satisfied with $\AM_{t}$ for the function $f$, which is equivalent to the validity of the inequality
\Eq{ac}{
f(tx+(1-t)y)\leq tf(x)\uplus (1-t)f(y).
}
If $\max\{f(x),f(y)\}=+\infty$, then the right hand side of \eq{ac} is equal to $+\infty$, thus, we can 
suppose that the right hand side is finite, that is $f(x)=h(x)$ and $f(y)=h(y)$. Now we have that $x\in \Q_0$ 
and $y\in\Q_0\cup\Q_1$. Then, using \eq{QQ}, it follows that 
$tx+(1-t)y\in \Q_0$. Therefore, applying the convexity of $h$, we get
\Eq{*}{
f(tx+(1-t)y)=h(tx+(1-t)y)\leq th(x)+(1-t)h(y)=tf(x)\uplus (1-t)f(y).
}
This proves that $f$ is upper $\AM_{t}$-convex for all $t\in\,]0,1[\,\cap\,\Q_1$.

To show that $f$ is not upper $\AM_{1-t}$-convex, let $y:=a\in\Q_1$ and let $x\in I\cap\Q_0$ be an arbitrary 
point. It follows from \eq{QQ} that the convex combination $(1-t)x+ty$ belongs to $\Q_1$ and 
it is also different from $a$. Therefore we have $f((1-t)x+ty)=+\infty$ and
$(1-t)f(x)\uplus tf(y)=(1-t)h(x)+th(y)\in\R$, which means that \eq{ac} cannot be satisfied.
\end{proof}

\Defn{Mf}{
For a function $f:I\to\RX$, define the following two classes of means:
\Eq{*}{
\underline{\M}_{f}&:=
\{M:I_\leq^2\to\R\mid M\text{ is a strict mean and }f\text{ is lower }M\text{-convex}\}, 
\\[1mm]
\overline{\M}_{f}&:=
\{M:I_\leq^2\to\R\mid M\text{ is a strict mean and }f\text{ is upper }M\text{-convex}\}.
}}

Note, that, due to the strictness of the means in the definition, the above sets can be also empty. The 
following proposition shows a certain algebraic closedness 
property of the classes $\underline{\M}_{f}$ and $\overline{\M}_{f}$.

\Prp{4p}{
For a function $f:I\to\RX$, the following statements hold:
\begin{enumerate}[(a)]\itemsep=0.2cm
\item\label{4p1} if $M,N_{1},N_{2}\in\underline{\M}_{f}$ (resp.\ $M,N_{1},N_{2}\in\overline{\M}_{f}$) and 
$N_{1}<N_{2}$ on the set $I_<^2$, then $M\circ(N_{1}, N_{2})\in\underline{\M}_{f}$ 
(resp.\ $M\circ(N_{1},N_{2})\in\overline{\M}_{f}$), and
\item\label{4p2} if $M, N\in\underline{\M}_{f}$ (resp.\ $M, N\in\overline{\M}_{f}$), then $M\circ(\min, N)$ 
and 
$M\circ(N, \max)$ also belong to $\underline{\M}_{f}$ (resp.\ to $\overline{\M}_{f}$).
\end{enumerate}
}
 
\begin{proof}
We verify the statements for the class $\overline{\M}_{f}$ only. The proof in the other case is completely 
analogous and also based on \lem{EFC}.

Let $(x,y)\in I_<^2$ be arbitrarily fixed, furthermore consider the points $p_{1}:=N_{1}(x,y)$ and 
$p_{2}:=N_{2}(x,y)$. (Obviously, under the conditions of (\ref{4p1}), it follows that $p_1<p_2$.) 
Using these notations, in view of \lem{EFC}, we need to show, that
\Eq{ee4}{
f(M(p_{1}, p_{2})) \leq
\frac{y-M(p_{1}, p_{2})}{y-x}f(x)\uplus \frac{M(p_{1}, p_{2})-x}{y-x}f(y),
}
holds. By applying the $M$- and then the $N_1$- and $N_2$-convexity of $f$, we have 
the 
following calculation:
\Eq{*}{
 f&(M(p_{1},p_{2}))\\ 
 &\leq\frac{p_{2}-M(p_{1},p_{2})}{p_2-p_1}f(p_{1})\uplus \frac{M(p_{1},p_{2})-p_{1}}{p_2-p_1}f(p_{2}) \\
 &=\frac{p_{2}-M(p_{1},p_{2})}{p_2-p_1}f(N_{1}(x,y))\uplus 
       \frac{M(p_{1},p_{2})-p_{1}}{p_2-p_1}f(N_{2}(x,y)) \\
 &\leq\frac{p_{2}-M(p_{1},p_{2})}{p_2-p_1}\left(\frac{y-p_{1}}{y-x}f(x)\uplus 
         \frac{p_{1}-x}{y-x}f(y)\right)
 \uplus \frac{M(p_{1},p_{2})-p_{1}}{p_2-p_1}\left(\frac{y-p_{2}}{y-x}f(x)\uplus 
         \frac{p_{2}-x}{y-x}f(y)\right) \\
 &=\frac{y-M(p_{1}, p_{2})}{y-x}f(x)\uplus \frac{M(p_{1}, p_{2})-x}{y-x}f(y).
   }
Thus the inequality \eq{ee4} is satisfied, which means the statement (a) is true. 

A completely similar calculation shows that the 
statement (b) is also valid.
\end{proof}

\Cor{c1}{
For a function $f:I\to\RX$, the classes
\Eq{*}{
\underline{\M}^{*}_f
&:=\{M\in\underline{\M}_{f}\mid M\text{ is separately continuous in both variables}\},\\[1mm]
\overline{\M}^{*}_f
&:=\{M\in\overline{\M}_{f}\mid M\text{ is separately continuous in both variables}\}
}
have no isolated points with respect to the pointwise convergence, namely for all $M\in\underline{\M}^{*}_f$ 
(resp.\ $M\in\overline{\M}^{*}_f$) there exist sequences of means $(L_{n}), (U_{n})\subseteq 
\underline{\M}^{*}_f$ (resp.\ $(L_{n}),(U_{n})\subseteq\overline{\M}^{*}_f$), such that $L_{n}<M<U_{n}$ for 
all $n\in\N$, furthermore $L_{n}\to M$ and $U_{n}\to M$ pointwise on $I_<^2$ as $n\to\infty$.
}
  
\begin{proof}
We prove the statement only for the class $\underline{\M}^{*}_f$.
 
Let $M\in\underline{\M}_{f}^{*}$ be an arbitrarily fixed mean. We show only that the sequence 
$(U_{n})$ exists, because the existence of $(L_{n})$ can be proved similarly.

Let $U_0=\max$ and, for $n\geq 1$, let $U_n:=M\circ(M, U_{n-1})$. In the first step we show
that the sequence $(U_{n})$ belongs to $\underline{\M}_{f}^{*}$. To see this, we prove, by induction, that 
$M<U_n<U_{n-1}$ for all $n\in\N$ on $I_<^2$. Let $(x,y)\in I_<^2$ be fixed. For $n=1$, 
using that $M$ is a strict mean, we get
\Eq{*}{
U_{1}(x,y)=M(M(x,y),U_0(x,y))=M(M(x,y), y)\in \,]M(x,y),y[\,=\,]M(x,y),U_0(x,y)[\,.
}
Assume that $M<U_n<U_{n-1}$ hold on $I_<^2$ for some $n\geq 2$. Using this assumption, for 
$n+1$, we obtain that
\Eq{*}{
U_{n+1}(x,y)=M(M(x,y),U_{n}(x,y))\in \,]M(x,y),U_n(x,y)[\,.
}
Hence $M(x,y)<U_{n+1}(x,y)<U_{n}(x,y)$ follows for all $(x,y)\in I_<^2$, which completes the proof of 
the induction. Thus, due to the \prp{4p}, it follows that $(U_{n})\subseteq\underline{\M}_{f}$. Moreover, 
by the definition, $U_{n}$ is a strict mean and separately continuous in both variables for all $n\in\N$, 
hence $(U_{n})\subseteq\underline{\M}_{f}^{*}$.

In the second step we show, that $U_{n}\downarrow M$ pointwise on $I_<^2$ as $n\to\infty$.
Let $(x,y)\in I_<^2$ be arbitrarily fixed again. Obviously, the sequence $(U_{n}(x,y))\subseteq\,]x,y[\,$ 
is convergent, because it is monotone decreasing and bounded from below by $M(x,y)$. Denote
$\lim_{n\to\infty}U_{n}(x,y)$ by $U^*(x,y)$ which, of course, cannot be smaller than $M(x,y)$. Upon taking 
the limit 
$n\to\infty$ in the identity
\Eq{*}{
U_{n}(x,y)=M(M(x,y),U_n(x,y)),
}
we get that
\Eq{*}{
U^*(x,y)=M(M(x,y),U^*(x,y)).
}
The inequality $M(x,y)<U^*(x,y)$ would contradict the strictness of $M$, therefore, $U^*(x,y)=M(x,y)$ must be 
valid.
\end{proof}

The following theorem is one of the main results of this paper. Roughly speaking, it states that the lower 
$M$-convexity property is inherited by the descendants.

\Thm{com}{
Let $f:I\to\RX$, $n\geq 2$ and $M_1,\dots, M_n\in\underline{\M}_f$ be continuous 
strict means. Then, for all $i\in\{1,\dots,n\}$, we have $\D_i(M_1,\dots, M_n)\subseteq\underline{\M}_f$.
}

\begin{proof}
Let $i\in\{1,\dots,n\}$ and $N\in\D_i(M_1,\dots, M_n)$ be arbitrarily fixed. We have already seen that, 
under our conditions, $N$ is a strict mean. If $(x,y)\in I_<^2$, then there exists $k\in\{1,\dots, n\}$ and 
$(\xi_1,\dots, \xi_n)\in\Phi_{(x,y)}$ such that $N(x,y)=\xi_k$, furthermore, with $\xi_0:=x$ and 
$\xi_{n+1}:=y$, we have
\Eq{*}{
M_j(\xi_{j-1}, \xi_{j+1})=\xi_j
\qquad (j\in\{1,\dots, n\}).
}
Using this and, for all $j\in\{1,\dots, n\}$, the lower $M_j$-convexity of the function $f$, we obtain
\Eq{*}{
0\leq \lfloor\xi_{j-1}, \xi_j, \xi_{j+1}; f\rfloor\qquad (j\in\{1,\dots, n\}).
}
Now, applying the Extended Chain Inequality, we get that
\Eq{*}{
0\leq \min_{1\,\leq\,j\,\leq\,n}\lfloor \xi_{j-1}, \xi_{j}, \xi_{j+1}; f\rfloor
  \leq \lfloor x, \xi_k, y; f\rfloor
  =\lfloor x, N(x,y), y; f\rfloor.
}
This means, by the definition, that $f$ is lower $N$-convex, that is $N\in\underline{\M}_{f}$.
\end{proof}

\Cor{com2}{
Let $f:I\to\RX$, $n\geq 2$, $s_{1},\dots, s_n\in\,]0,1[\,$, and let $h:I\to\R$ be a continuous, 
strictly increasing function. Assume that $\MM_{s_ih,\,(1-s_i)h}\in\underline{\M}_f$ for all $i\in\{1,\dots, 
n\}$.
Then, for all $i\in\{1,\dots, n\}$, the Matkowski mean $\MM_{\sigma_ih,\,(1-\sigma_i)h}$ also
belongs to $\underline{\M}_f$, where
\Eq{si2}{
 \sigma_i:=\bigg(\sum_{j=i}^n\prod_{k=1}^{j}\frac{s_k}{1-s_k}\bigg)
   \bigg(\sum_{j=0}^n\prod_{k=1}^{j}\frac{s_k}{1-s_k}\bigg)^{-1} \qquad(i\in\{1,\dots,n\}).
}}

\begin{proof}
For $(x,y)\in I_<^2$, define the mapping $\varphi_{(x,y)}:[x,y]_\leq^{n}\to\R^{n}$ as in \thm{cmat}.
In view of this theorem, it follows that, for all $(x,y)\in I_<^2$, the fixed point set 
$\Phi_{(x,y)}$ is the singleton $\{(\xi_{1},\dots, \xi_{n})\}$, where 
$\xi_i=\MM_{\sigma_ih,\,(1-\sigma_i)h}(x,y)$. Thus, for $i\in\{1,\dots, n\}$, 
the function $\MM_{\sigma_ih,\,(1-\sigma_i)h}$ is the $i^{\mbox{\tiny\rm th}}$ descendant of the $n$-tuple of 
means $(\MM_{s_1h,\,(1-s_1)h},\dots,\MM_{s_nh,\,(1-s_n)h})$. Therefore, due to \thm{com}, we obtain that 
$\MM_{\sigma_ih,\,(1-\sigma_i)h}\in\underline{\M}_f$ for all $i\in\{1,\dots, n\}$.
\end{proof}

\Cor{wam1}{
Let $n\geq 2$, $p, q, h_{1},\dots, h_{n-1}:I\to\R$ be continuous, strictly increasing 
functions and $f:I\to\RX$. Set further $h_0:=h_n:=0$ and assume that there exists 
$j\in\{1,\dots, n\}$ such that
\Eq{*}{
\big\{\MM_{p+h_{i-1},\,h_i}\mid 1\leq i\leq j-1\big\}\cup
   \big\{\MM_{p+h_{j-1},\,q+h_j}\big\}\cup
   \big\{\MM_{h_{i-1},\,q+h_{i}}\mid j+1\leq i\leq n\big\}
\subseteq\underline{\M}_{f}.
}
Then $N_{1},\dots, N_{n}\in\underline{\M}_{f}$, where, for all $(x,y)\in I_\leq^2$,
\Eq{*}{
N_{j}(x,y)=\MM_{p,\,q}(x,y) \qquad\text{and}\qquad
N_i(x,y)=
\begin{cases}
 \MM_{p,\,h_i}(x, N_{i+1}(x,y)) & \text{if } i\in\{1,\dots, j-1\},\\[2mm]
 \MM_{h_{i-1},\,q}(N_{i-1}(x,y), y) & \text{if } i\in\{j+1,\dots, n\}.
\end{cases}}}

\begin{proof}
The method of the proof is same as that of \cor{com2}. For $(x,y)\in I_<^2$, 
define the mapping $\varphi_{(x,y)}$ as in \eq{phid} by the using the means 
$M_{1},\dots, M_{n}$, where
\Eq{*}{
M_i:=
\begin{cases}
\MM_{p+h_{i-1},\,h_i} & \text{if } i\in\{1,\dots, j-1\}, \\[2mm]
\MM_{p+h_{i-1},\,h_{i}+q} & \text{if }i=j, \\[2mm]
\MM_{h_{i-1},\,h_{i}+q} & \text{if } i\in\{j+1,\dots, n\}.
\end{cases}
}
Due to \thm{rmat}, it follows that, for all $(x,y)\in I_<^2$, the fixed point set 
$\Phi_{(x,y)}$ is the singleton $\{(\xi_{1},\dots, \xi_{n})\}$, where we have
\Eq{*}{
\xi_j:=\MM_{p,\,q}(x,y)
\qquad\text{and}\qquad
\xi_i:=
 \begin{cases}
  \MM_{p,\,h_{i}}(x,\xi_{i+1}) &\text{if }i\in\{1,\dots, j-1\}, \\[2mm]
  \MM_{h_{i-1},\,q}(\xi_{i-1},y) &\text{if }i\in\{j+1,\dots, n\}.
 \end{cases}
}
Thus, for $i\in\{1,\dots, n\}$, the function $N_i:I_\leq^2\to\R,\,N_i(x,y):=\xi_i$
is the $i^{\mbox{\tiny\rm th}}$ descendant of the $n$-tuple of means
$(\MM_{p+h_{i-1},\,h_i},\dots, \MM_{p+h_{j-1},\,h_{j}+q},\dots, \MM_{h_{i-1},\,h_{i}+q})$.
Hence, by \thm{com}, it follows that $N_{i}\in\underline{\M}_{f}$ for all $i\in\{1,\dots, n\}$.
\end{proof}

\section{$\AM_t$-convexity of extended real valued functions}

In this section we investigate a special subclass of $\underline{\M}_f$ and $\overline{\M}_f$, respectively.
For an extended real valued function $f:I\to\RX$ consider the sets $\underline{\AC}_f$ and $\overline{\AC}_f$
defined by
\Eq{*}{
\underline{\AC}_f&:=
\{0<t<1\mid\text{ for all $(x,y)\in I_<^2$ : }\lfloor x, \AM_t(x,y), y;f\rfloor\geq 0\},\\
}
and
\Eq{*}{
\overline{\AC}_f&:=
\{0<t<1\mid\text{ for all $(x,y)\in I_<^2$ : }\lceil x, \AM_t(x,y), y;f\rceil\geq 0\}.
}
If $f$ is real-valued, then clearly these two sets are the same, therefore, we will simply denote them by 
$\AC_f$. Note that, by the definitions, both sets can also be empty. On the other hand, these sets can be
easily identified with the subclass of weighted arithmetic means in $\underline{\M}_f$ and 
$\overline{\M}_f$ respectively, more precisely we have the following identifications
\Eq{*}{
t\in\underline{\AC}_f \quad\Longleftrightarrow\quad \AM_t|_{I_\leq^2}\in\underline{\M}_{f}
\qquad\mbox{and}\qquad
t\in\overline{\AC}_f \quad\Longleftrightarrow\quad \AM_t|_{I_\leq^2}\in\overline{\M}_{f}.
}
The motivation for our investigations is a well known result, which is due to N. Kuhn \cite{Kuh84}, 
and which is about the structure of the set of parameters for which a given real valued function is 
convex. The theorem says that if $f:I\to\R$ is an arbitrary function and
\Eq{*}{
\AC^\circ_f
  :=\AC_f\cap(1-\AC_f)=\{0<t<1\mid f\text{ is simultaneously $\AM_t$-convex and $\AM_{1-t}$-convex on }I\},
}
then we have that either $\AC^\circ_f=\emptyset$ or $\AC^\circ_f=K\cap\,]0,1[\,$, where $K$ is a subfield of 
$\R$. Moreover, the reverse of this statement is also valid: if $K\subseteq\R$ is a given subfield, then 
there exists a function $f:I\to\R$ such that $\AC^\circ_f$ equals to the intersection $K\cap\,]0,1[\,$.

The following results are about such algebraical closedness properties of the sets $\underline{\AC}_f$
and $\overline{\AC}_f$.

\Thm{tac1}{
Given a function $f:I\to\RX$, the following statements hold for $S\in\{\underline{\AC}_f,\overline{\AC}_f\}$:
\begin{enumerate}\itemsep=2mm
 \item\label{cac1A} if $t, s_1, s_2\in S$ with $s_1<s_2$, then $ts_{2}+(1-t)s_1\in S$,
 \item\label{cac1B} if $t,s\in S$, then $ts$ and $1-(1-t)(1-s)$ belong to $S$, and
 \item\label{cac1a} $S$ is dense in the open unit interval, provided that it is not empty.
\end{enumerate}}

\begin{proof}
We verify only the statements about $S=\underline{\AC}_f$. The proof for $S=\overline{\AC}_f$ is analogous.

Let $t, s_1, s_2\in \underline{\AC}_f$ with $s_1<s_2$. Then the means $\AM_t, \AM_{s_1}$ and $\AM_{s_2}$ 
belong to 
$\underline{\M}_f$ and, because of $s_1<s_2$, we have $\AM_{s_2}<\AM_{s_1}$ on $I_<^2$. Using \prp{4p} for 
$M:=\AM_t$, $N_1:=\AM_{s_2}$ and $N_2:=\AM_{s_1}$, we obtain that $\AM_t\circ(\AM_{s_2},
\AM_{s_1})\in\underline{\M}_{f}$. 
On the other hand, for $(x,y)\in I_<^2$, we have
\Eq{*}{
\AM_t\circ(\AM_{s_2}, \AM_{s_1})(x,y)&=\AM_t(\AM_{s_2}(x,y), \AM_{s_1}(x,y))=\AM_t(s_2x+(1-s_2)y, 
s_1x+(1-s_1)y) \\
&=(ts_2+(1-t)s_{1})x+(1-(ts_2+(1-t)s_{1}))y=\AM_{ts_2+(1-t)s_1}(x,y),
}
consequently $ts_2+(1-t)s_1\in\underline{\AC}_{f}$, which proves (\ref{cac1A}).

To prove (\ref{cac1B}), observe that, under our notation, $\min=\AM_1$ and $\max=\AM_0$ on $I_\leq^2$. 
Thus, according to the second statement of \prp{4p}, the means $\AM_t\circ(\AM_1, \AM_s)$ and 
$\AM_t\circ(\AM_s, \AM_0)$ belong to $\underline{\M}_{f}$. Then the same calculation yields that
$1-(1-t)(1-s)$ and $ts$ belong to $\underline{\AC}_f$, respectively.

To verify (\ref{cac1a}) assume that $\underline{\AC}_{f}$ is not empty and
indirectly suppose that $\underline{\AC}_{f}$ is not dense in $\,]0,1[\,$, that is there exist
$\alpha<\beta$ in $[0,1]$ such that $\underline{\AC}_{f}\,\cap\,]\alpha,\beta[\,$ is empty. 
We may assume that the interval $\,]\alpha,\beta[\,$ is maximal, or equivalently, for all
$\varepsilon>0$, the intersection $\underline{\AC}_{f}\,\cap\,]\alpha-\varepsilon,\beta+\varepsilon[\,$
is not empty. Observe that, due to the second assertion of the theorem, it easily follows that $0<\alpha$ and 
$\beta<1$. Indeed, if $t\in\underline{\AC}_f$ is arbitrary, then, due to the fact that
$\underline{\AC}_f$ is closed under the multiplication, for all $k\in\N$, the value $t^{k}$ belongs again to 
$\underline{\AC}_f$. Thus any open neighborhood of zero contains an element from $\underline{\AC}_f$, which
means $0<\alpha$. Similarly, using the closedness of $\underline{\AC}_f$ under the operation 
$(t,s)\longmapsto 1-(1-t)(1-s)$, we get that $\beta<1$. Thus we obtained that
$\,[\alpha, \beta]\,\subseteq\,]0, 1[\,$. Now, let $t\in\underline{\AC}_f$ be arbitrarily fixed and 
$(r_n), (s_n)\subseteq\underline{\AC}_f$ be sequences such that $r_n\nearrow\alpha$ and $s_n\searrow\beta$ 
as $n\to\infty$. Then, in view of the first assertion of the theorem, $ts_n+(1-t)r_n\in\underline{\AC}_{f}$ 
for all $n\in\N$ and $ts_n+(1-t)r_n\to t\beta+(1-t)\alpha\in\,]\alpha,\beta[$ as $n\to\infty$. Therefore, for 
sufficiently large $n$, we get that $ts_n+(1-t)r_n\in ]\alpha,\beta[$, which contradicts the emptiness of 
$\underline{\AC}_{f}\,\cap\,]\alpha,\beta[\,$ and hence $\underline{\AC}_f$ must be dense in $\,]0, 1[\,$.
\end{proof}

\Rem{NE}{The result stated in \thm{tac1} is not analogous to that of Kuhn \cite{Kuh84}. In general, the set 
$\overline{\AC}_{f}$ is not of the form $]0,1[\cap K$, where $K$ is a subfield of $\R$. To see this, it is sufficient 
to construct a function $f:I\to\RX$ such that the set $\overline{\AC}_{f}$ is not closed under the addition of their 
elements. 

Indeed, let $f:I\to\R\cup\{+\infty\}$ be the function defined in \prp{x1}. For arbitrarily fixed parameters
$s,t\in\,]0, 1[\,\cap\,\Q_1\subseteq\overline{\AC}_f$ with $s+t<1$, in view of \eq{QQ}, the sum $s+t$ belongs to $\Q_0$.
To prove that $s+t\not\in \overline{\AC}_f$, we construct $x<y$ in $I$ such that
\Eq{spt}{
f((s+t)x+(1-(s+t))y)>(s+t)f(x)+(1-(s+t))f(y).
}
Let $x\in I\cap\Q_0$ be arbitrarily fixed and set $y:=a$. Then, using again \eq{QQ}, the convex combination
$u:=(s+t)x+(1-(s+t))y$ belongs to $I\cap\Q_1$ and it is also different from $a$. Consequently $f(u)=+\infty$, on the 
other hand
\Eq{*}{
(s+t)f(x)+(1-(s+t))f(y)=(s+t)h(x)+(1-(s+t))h(y)\in\R,
}
thus \eq{spt} is satisfied.}

\Cor{com3}{
Let $I\subseteq\R$ be an interval, $f:I\to\RX$, $n\geq 2$ and $s_{1},\dots, s_n\in\underline{\AC}_f$.
Then $\sigma_i\in\underline{\AC}_f$ for all $i\in\{1,\dots, n\}$, where
\Eq{si3}{
 \sigma_i:=\bigg(\sum_{j=i}^n\prod_{k=1}^{j}\frac{s_k}{1-s_k}\bigg)
   \bigg(\sum_{j=0}^n\prod_{k=1}^{j}\frac{s_k}{1-s_k}\bigg)^{-1}.
}}

\begin{proof}
Apply \cor{com2} under $h:=\id$.
\end{proof}

\Cor{com4}{
For a function $f:I\to\RX$ the following statements hold:
\begin{enumerate}\itemsep=2mm
\item\label{c4a} if $1/2\in\underline{\AC}_f$ then $r\in\underline{\AC}_f$ for all $r\in\Q\,\cap\,]0,1[\,$,
\item\label{c4b} if $\ell/m\in\underline{\AC}_f$ for some $\ell,m\in\N$ with $\ell<m$ and $\ell\neq m/2$, then, for 
all $n\geq 2$ and for all $i\in\{1,\dots, n\}$, the fraction
\Eq{*}{
r_i:=
  \frac{\ell^{n+1}-\ell^{i}(m-\ell)^{n+1-i}}
	     {\ell^{n+1}-(m-\ell)^{n+1}}
}
belongs to $\underline{\AC}_f$.
\end{enumerate}}

\begin{proof}
To prove (\ref{c4a}), assume that $1/2\in\underline{\AC}_f$ and let $p,q\in\N$ be arbitrarily
fixed numbers such that $q>1$ and $p<q$. For $q=2$, the statement (\ref{c4a}) is trivial, thus we may
assume that $q>2$. Now set $n:=q-1$ and $i_0:=q-p$. Then $n\geq 2$ and $i_0\in\{1,\dots, n\}$, 
thus, using \cor{com3} for $s_1:=\dots=s_n:=1/2$, we get that
\Eq{*}{
\sigma_{i_0}=\frac{n-i_0+1}{n+1}=\frac{q-1-(q-p)+1}{q-1+1}=\frac{p}{q}.
}
This means that $\Q\,\cap\,]0,1[\,\subseteq\underline{\AC}_f$.

To prove (\ref{c4b}), assume that $\ell/m\in\underline{\AC}_f$ for some $\ell,m\in\N$, where 
$\ell<m$ and $2\ell\neq m$. Let further $n\geq 2$ be arbitrarily fixed and set $s_1:=\dots=s_n:=\ell/m$.
Then a simple calculation shows that $\sigma_i=r_i$ for all $i\in\{1,\dots, n\}$. Due to
\cor{com3}, we get that $r_i\in\underline{\AC}_f$ for all $i\in\{1,\dots, n\}$.
\end{proof}


\begin{thebibliography}{1}

\bibitem{DarPal87}
Z.~Daróczy and Zs. Páles, \emph{{Convexity with given infinite weight
  sequences}}, Stochastica \textbf{11} (1987), no.~1, 5–12. \MR{90c:39020}

\bibitem{GilPal08}
A.~Gilányi and Zs. Páles, \emph{{On convex functions of higher order}}, Math.
  Inequal. Appl. \textbf{11} (2008), no.~2, 271–282.

\bibitem{Kuc85}
M.~Kuczma, \emph{{An {I}ntroduction to the {T}heory of {F}unctional {E}quations
  and {I}nequalities}}, {Prace Naukowe Uniwersytetu Śląskiego w Katowicach},
  vol. 489, Państwowe Wydawnictwo Naukowe — Uniwersytet Śląski,
  Warszawa–Kraków–Katowice, 1985, 2nd edn. (ed. by A. Gilányi),
  Birkhäuser, Basel, 2009. \MR{0788497 (86i:39008), MR 2467621}

\bibitem{Kuh84}
N.~Kuhn, \emph{{A note on $t$-convex functions}}, {General Inequalities, 4
  (Oberwolfach, 1983)} (W.~Walter, ed.), {International Series of Numerical
  Mathematics}, vol.~71, Birkhäuser, Basel, 1984, p.~269–276. \MR{87b:26015}

\bibitem{LewOlb14}
M.~Lewicki and A.~Olbryś, \emph{{On non-symmetric {$t$}-convex functions}},
  Math. Inequal. Appl. \textbf{17} (2014), no.~1, 95–100. \MR{3220977}

\bibitem{Mat10b}
J.~Matkowski, \emph{{Generalized weighted and quasi-arithmetic means}},
  Aequationes Math. \textbf{79} (2010), no.~3, 203–212. \MR{2665531}

\bibitem{NikPal03}
K.~Nikodem and Zs. Páles, \emph{{On $t$-convex functions}}, Real Anal.
  Exchange \textbf{29} (2003), no.~1, 219–228. \MR{2005f:26030}

\end{thebibliography}

\def\MR#1{}

\end{document}